\documentclass{amsart}[12pt,a4paper]

\usepackage[colorlinks=true,          
            linkcolor=blue,
            citecolor=red,
            urlcolor=blue]{hyperref}
\usepackage{amssymb,amsfonts,amsmath,hyperref,mathrsfs,latexsym,stmaryrd,graphicx,manfnt,enumitem,caption}
\DeclareMathAlphabet{\mathpzc}{OT1}{pzc}{m}{it}

\DeclareMathAlphabet{\mathpzc}{OT1}{pzc}{m}{it}

\makeatletter
\renewcommand*\env@matrix[1][c]{\hskip -\arraycolsep
  \let\@ifnextchar\new@ifnextchar
  \array{*\c@MaxMatrixCols #1}}
\makeatother

 \newtheorem{thm}{Theorem}[section]
  \theoremstyle{plain}

 \newtheorem{Prop}{Proposition}[section]
  \theoremstyle{definition}

\newcommand {\CC}{\mathbb{C}}

\newcommand {\RR}{\mathbb{R}}

\newcommand{\ZZ}{\mathbb{Z}}

\newcommand {\bal}{\boldsymbol{\alpha}}
\newcommand {\bbe}{\boldsymbol{\beta}}

\newcommand {\bla}{\boldsymbol {\lambda}}

\newcommand {\bv}{{\bf v}}

\newcommand {\bt}{{\bf t}}
\newcommand {\bu}{{\bf u}}
\newcommand {\bw}{{\bf w}}

\newcommand {\bI}{{\bf I}}

\newcommand {\bdelta}{\boldsymbol\delta}

\newcommand{\cB}{\mathcal{B}}

\newcommand{\cL}{\mathcal{L}}

\newcommand {\cO}{\mathcal{O}}
\newcommand{\cP}{\mathcal{P}}

\newcommand{\cT}{\mathcal{T}}

\newcommand{\cV}{\mathcal{V}}

\newcommand{\cX}{\mathcal{X}}
\newcommand{\cY}{\mathcal{Y}}

\newcommand{\scB}{\mathscr{B}}

\newcommand{\fg}{\mathfrak{g}}

\newcommand{\fl}{\mathfrak{l}}

\newcommand{\fo}{\mathfrak{o}}

\newcommand{\fs}{\mathfrak{s}}

\newcommand{\ft}{\mathfrak{t}}

\newcommand{\fD}{\mathfrak{D}}

\newcommand{\sym}{\textrm{Sym}}

\newcommand{\mhom}{\textrm{Hom}}

\newcommand{\pr}{\textrm{pr}}

\newcommand{\tot}{\textrm{tot }}

\newcommand{\ctimes}{\otimes_\CC}

\newcommand {\io}{\iota}

\newcommand{\Prym}{{\bf Prym}}

\newcommand{\Pic}{{\bf Pic}}

\newcommand{\Higgs}{{\bf Higgs}}

\newcommand{\hookr}{\hookrightarrow}

\include{biblio}

\usepackage{latexsym}
\usepackage{mathrsfs}
\input xy
\xyoption{all}
\usepackage{tikz}
\usepackage{tikz-3dplot}
\usepackage{tkz-fct}
\usepackage{tikz-cd}
\usepackage{tkz-euclide}
\usepackage[framemethod=tikz]{mdframed}

\title{ Donagi--Markman cubics for Hitchin systems of type $A_2$, $B_2$, $G_2$ }
\author[U. Bruzzo and P. Dalakov]{Ugo Bruzzo$^{\P\sharp\star}$ and Peter Dalakov$^{\dag\ddag}$}
\address{\small $^\P$ SISSA (Scuola Internazionale Superiore di Studi Avanzati)\\ 
Via Bonomea 265, 34136 Trieste, Italia\\
$^\sharp$ INFN (Istituto Nazionale di Fisica Nucleare), Sezione di Trieste \\
$^\star$ IGAP (Institute for Geometry and Physics), Trieste, Italy\\
$^\dag$ American University in Bulgaria, 2700 Blagoevgrad, G.Izmirliev Sq.1, Bulgaria\\
$^\ddag$ Institute of Mathematics and Informatics, Bulgarian Academy of Sciences, Sofia, Bulgaria}

\date{\today}

\begin{document}
	\begin{abstract}
		We obtain explicit formulae for  the Donagi--Markman (Bryant--Griffiths, Yukawa) cubic for Hitchin
		systems of type $A_2$, $B_2$ and $G_2$. This is achieved by evaluating the quadratic residues in
		the Balduzzi--Pantev formula, using a previous result of ours. For $G_2$ we also recover earlier results   of Hitchin.
	\end{abstract}

\maketitle
\tableofcontents

	\section{Introduction}

Let $X$ be a smooth connected Riemann surface of genus $g\geq 2$.
We  write $\cB_{\fg_2}$ or just $\cB$  for the  \emph{$G_2$ Hitchin base}, $\cB_{\fg_2}:=H^0(X,K_X^2)\oplus H^0(X,K_X^6)$.
We shall think of $\cB_{\fg_2}$ as the base of the universal family 
$p:\cX\to \cB_{\fg_2}$ 
of ($K_X$-valued)
$G_2$-cameral covers of $X$, see section \ref{section_cam_curve}. We also write $p_b:\widetilde{X}_b\to X$ for the cover, corresponding to $b\in \cB_{\fg_2}$.
We are going to denote by
$\scB\subseteq \cB$ the Zariski open subset of smooth cameral curves with generic ramification.

There is a family of abelian varieties, $\Prym_{\cX/\scB}\to \scB$, whose fibres are generalized (cameral) Prym varieties. Explicitly, the
(generalized)
Prym  variety, associated to the cover
 $p_b:\widetilde{X}_b\to X$ is $\left(\Lambda\otimes_\ZZ\Pic_{\widetilde{X}_b}\right)^W= H^1(X,\cT_b)$,  where $\Lambda$ is the cocharacter lattice
 of $G_2$, $W=D_6$ is its Weyl group and $\cT_b=(p_b)_\ast(\Lambda\otimes \cO_{\widetilde{X}_b}^\times)^W$.
 The family of cameral Pryms is in fact an algebraic completely integrable Hamiltonian system -- so, in particular, its fibres are (complex) Lagrangian,
 see e.g., \cite[\S 5.3]{hurtubise_markman_rk2} for a general construction.

 The base $\cB_{\fg_2}$ is also the base of the ($G_2$) Hitchin  integrable system $h: \Higgs_{G_2}^{ss}\to \cB_{\fg_2}$,
 whose total space is the moduli space of semi-stable $K_X$-valued (principal) $G_2$ Higgs bundles on $X$, see
 \cite{hitchin_SB}, \cite{katz_pan_G2}, \cite{hitchin_G2}.
 The two integrable systems are related by
 ``abelianization''.   
 In particular, a Hitchin fibre $h^{-1}(b)$, $b\in \scB$, is a torsor for $H^1(X,\cT_b)$. Some of the  references for this
 much studied relation are
 \cite{katz_pan_G2},\cite{scognamillo_elem}, \cite{faltings},\cite{donagi_spectral_covers}, \cite{don-gaits},\cite{hitchin_G2}.
 While this relation is motivationally important for us, we are not going to review it, as the two
 integrable systems have the same infinitesimal period map.

 The infinitesimal period map for a family of abelian varieties having the structure of a Lagrangian fibration satisfies
 the ``cubic condition'' of Donagi and Markman (\cite{donagi_markman}, \cite{donagi_markman_cubic}) and determines a section
   $c\in H^0(\scB, \sym^3T^\vee_\scB)$. We recall  the main ideas and fix notation in section \ref{section_periods}.
  For the case of Hitchin Pryms, the Balduzzi--Pantev formula (\cite{balduzzi}), combined with with \cite[Proposition 2.8] {hertling_hoev_posthum} gives us that
	  $c$, when evaluated on a triple of tangent vectors $\bu,\bv,\bw\in T_b\scB=\cB$, for $b\in\scB$, is given by

	  \begin{equation}\label{BP_ver1}
	      c_b(\bu,\bv,\bw)=\frac{1}{2}\sum_{x\in \textrm{Ram }p_b}\textrm{Res}^2_x
	  \left(p_b^\ast \frac{\cL_{\bu} \mathfrak{D}(b)}{\mathfrak{D}(b)}\left(\nabla^{GM}_{\bv}\bla_{SW}\right)_b\cup \left(\nabla^{GM}_{\bw}\bla_{SW}\right)_b \right),
	  \end{equation}
	  see also \cite[eq.~(8)]{bruzzo_dalakov_sw}.  Here $\bla_{SW}$ is the Seiberg--Witten differential and $\nabla^{GM}$
	     the Gauss--Manin connection, see section \ref{section_periods}.
	     We review the discriminant $\fD$ in section \ref{section_reminder_G2}, but recall that it factors as
	     $\fD=\fD_1\fD_2$, with the two factors corresponding to the products of roots of different lengths.
	    The formula (\ref{BP_ver1}) holds, more generally,  for $K_X$-valued $W$-cameral covers,
	  where $W$ is the Weyl group of a 
	  reductive Lie group $G$ (with chosen Borel and Cartan subgroups). It should be noted that in such a generality,
	  the associated Prym varieties and Higgs bundle moduli spaces depend in an essential way on 
	   the group  $G$,  see \cite{don-gaits}, \cite{scognamillo_elem}. On the other hand, the Hitchin base 
	   and the Donagi--Markman cubic  depend  on $G$ only through $\fg=\textrm{Lie }G$.
	  
	  In this note, we use
	  \cite[Theorem A]{bruzzo_dalakov_sw} and a careful analysis of the 12-fold covers $p_b:\widetilde{X}_b\to X$ to 
	     calculate $\nabla^{GM}_\bv\bla_{SW}$ at ramification points and evaluate  the quadratic residues
	     in  formula (\ref{BP_ver1})  explicitly.

	  \begin{thm}\label{thm1}
			At $b=(b_1,b_2)\in \scB_{\fg_2}$ the value of the $G_2$ Donagi--Markman cubic is
			\begin{multline}\label{cubic_G2_fancy}
		      	c_b= 36\sum_{z_0\in V(\fD_1(b))}\textrm{ev}_{z_0}\frac{(d\fD_1)_b\otimes (d\fD_1)_b\otimes (d\fD_1)_b}{b_1^2 j^1(\fD_1(b))^2} \\
		      	+ \frac{4}{3}\sum_{z_0\in V(\fD_2(b))}\textrm{ev}_{z_0}\frac{(d\fD_2)_b\otimes (d\fD_2)_b\otimes (d\fD_2)_b}{b_1^2 j^1(\fD_2(b))^2}.
		      \end{multline}
			
			Explicitly, for a given triple of
			 tangent vectors
			\[
			    \bu=(u_1,u_2), \bv=(v_1,v_2), \bw=(w_1,w_2)\in  T_b\scB=H^0(X,K^2)\oplus H^0(X,K^6)
			\]
			one has 
			\begin{multline}
				c_b(\bu,\bv,\bw)= 36\sum_{z_0\in V(b_2)}\frac{u_2 v_2 w_2}{b_1^2 j^1 (b_2)^2}(z_0)+\\
				36 \sum_{z_0\in V(4b_1^3-27b_2)} \frac{ \left(-u_2+\frac{4}{9}b_1^2u_1\right)\left(-v_2+\frac{4}{9}b_1^2v_1\right)\left(-w_2+\frac{4}{9}b_1^2w_1\right)}
				{b_1^2 j^1(-b_2+\frac{4}{27}b_1^3)^2}(z_0).
			\end{multline}
			
		\end{thm}

		We recall that the Atiyah (first jet) sequence for a vector bundle $E$ on  the curve $X$ is 
		\begin{equation}\label{Atiyah_sqs}
			\xymatrix@1{0\ar[r] & K_X\otimes E\ar[r] & J^1(E)\ar[r]& E\ar[r]& 0},
		\end{equation}
		and for $E=K_X^n$ this is
		\begin{equation}\label{jet_canonical}
			\xymatrix@1{0\ar[r] & K_X^{n+1}\ar[r] & J^1(K_X^n)\ar[r]& K^n_X\ar[r]& 0}.
		\end{equation}
		Hence the ratios in both summands of  the formula make sense, with numerator and denominator being at $z_0$ 
		in the fibre  of the same line bundle ($K_X^{18}$ and $K_X^{24}$,
		respectively).
		In  section \ref{section_pf} we provide  some equivalent formulations and  local expressions, see e.g. equations (\ref{residues_final_global}), 
		(\ref{residues_final_local}).

		While a calculation of the $G_2$ Donagi--Markman cubic appears in \cite{hitchin_G2}, we provide a different  global expression for it
		 and a different
		 proof.

		There is an additional symmetry of the $G_2$ cubic, related to the Langlands involution of $\cB_{\fg_2}$, as we now briefly recall.

		It is known that  Hitchin systems have a certain duality property, see \cite{don-pan},  \cite{hitchin_G2},
		\cite{tad-haus} and the physics papers  \cite{KW}, \cite{argyres_kapustin_seiberg}. 
		Let $G$ and $^LG$ be a pair of Langlands-dual complex reductive Lie groups and let $\cB$ and $^L\cB$ be the respective Hitchin bases.
		By \cite[Theorem A]{don-pan}, there is an isomorphism $\tau:\cB\simeq ^L\cB$, 
		which   preserves the discriminants and lifts to an isomorphism between
		the universal cameral covers.
		Away from the discriminant loci, the Hitchin fibres over pairs of points, matched by this isomorphism, are dual. The duality is induced by an isomorphism of
		polarised abelian varieties, and is in fact the restriction of a global duality.
		For the case of Langlands self-dual Lie groups, the isomorphism between the respective bases can be considered 
		as an involution of $\cB$
		(the two Hitchin bases being naturally isomorphic). If $G$ is simple,   this automorphism
		is different from the identity only in the non simply-laced case, i.e., in the cases of $G_2$ and $F_4$ (types
		$B_n$ and $C_n$ are not self-dual, they are Langlands dual to each other).
		We now focus on the case of $G_2$, and denote by $\tau: \cB_{\fg_2}\to \cB_{\fg_2}$   the Langlands  involution of $\cB_{\fg_2}$, 
		as well for its restriction to $\scB$. For the explicit formula (with our conventions),
		see section \ref{section_reminder_G2}. We write $\tau_i:\scB\to H^0(X, K_X^{d_i})$, $i=1,2$
		for the components $\textrm{pr}_i\circ \tau$ of $\tau$, see equations (\ref{tau_base}) and (\ref{tau_components})
		for explicit formulae. We notice that with our conventions $\fD_2=27 \tau_2$.

		\begin{thm}\label{thm2}
			In the above notation, we have
			\begin{multline}
				c_b(\bu,\bv,\bw)= 36\sum_{z_0\in V(b_2)}\frac{u_2 v_2 w_2}{b_1^2 j^1 (b_2)^2}(z_0)+\\
				36 \sum_{z_0\in V(\tau_2(b))} \frac{ d\tau_{2 b}(\bu)d\tau_{2 b}(\bv) d\tau_{2 b}(\bw)}{b_1^2 j^1(\tau_2(b))^2}(z_0)
			\end{multline}
			and hence
			\begin{equation}
      	\tau^\ast c=c\in H^0(\scB,\sym^3 T^\vee_{\scB}).
	  \end{equation}
		\end{thm}
		
		We thus establish that  the $G_2$ Donagi--Markman cubic is invariant under the Langlands involution,
		 a result appearing in \S 6.5 of \cite{hitchin_G2}.

		Our method, based on \cite{bruzzo_dalakov_sw}, is not limited to $G_2$ in any way. 
                In the two Appendices we consider the other two simple Lie algebras of rank two and obtain formulae for the
                Donagi--Markman cubics.
                Since the structure of the argument is the same as in $G_2$ case, we only indicate the main steps 
                of the derivations without providing exhaustive details.
                
		First,  in  Appendix \ref{appendix_A2} we consider  $A_2=\fs\fl_3(\CC)$ Hitchin systems.
		Here  the Hitchin base is $\cB_{\fs\fl_3}=H^0(X, K_X^2)\oplus H^0(X,K_X^3)$, and we write again $\scB_{\fs\fl_3}$ for the Zariski open subset
		of generic $A_2$ cameral covers.
		
		 \begin{thm}\label{thm3}
			At  $b=(b_1,b_2)\in \scB_{\fs\fl_3}$, 
			the value of the  $A_2$ Donagi--Markman cubic is
			\begin{equation}\label{cubic_A2_fancy}
		      	c_b= \frac{9}{2}\sum_{z_0\in V(\fD(b))}\textrm{ev}_{z_0}\frac{d\fD_b\otimes d\fD_b\otimes d\fD_b}{b_1^2 j^1(\fD(b))^2}.
		      \end{equation}
			
			Explicitly, given a triple of tangent vectors
			\[
			    \bu=(u_1,u_2), \bv=(v_1,v_2), \bw=(w_1,w_2)\in  T_b\scB_{\fs\fl_3}=H^0(X,K_X^2)\oplus H^0(X,K_X^3)
			\]
			one has
			\begin{multline}
				c_b(\bu,\bv,\bw)= \\  -\frac{1}{6} \sum_{z_0\in V(4b_1^3-b^2_2)}\frac{(12b_1^2u_1-2b_2 u_2)(12b_1^2 v_1-2b_2v_2)(12b_1^2 w_1- 2b_2 w_2)}{b_1^2 j^1 (4b_1^3-b_2^2)^2}(z_0).
			\end{multline}

		\end{thm}

		Finally, in Appendix \ref{appendix_B2} we discuss the case of  $B_2=\fs\fo_5(\CC)$ Hitchin systems. Here the
		Hitchin base is $\cB_{\fs\fo_5}= H^0(X,K_X^2)\oplus H^0(X,K_X^4)$ and $\scB_{\fs\fo_5}$ is the Zariski open subset of
		generic $B_2$ cameral covers.
		
		\begin{thm}\label{thm4}
			At  $b=(b_1,b_2)\in \scB_{\fs\fo_5}$, 
			the value of the  $B_2$ Donagi--Markman cubic is
			\begin{multline}\label{cubic_B2_fancy}
			c_b= 12\sum_{z_0\in V(\fD_1(b))}\textrm{ev}_{z_0}\frac{(d\fD_1)_b\otimes (d\fD_1)_b\otimes (d\fD_1)_b}{b_1 j^1(\fD_1(b))^2} \\
		      	+ 3\sum_{z_0\in V(\fD_2(b))}\textrm{ev}_{z_0}\frac{(d\fD_2)_b\otimes (d\fD_2)_b\otimes (d\fD_2)_b}{b_1 j^1(\fD_2(b))^2}.
		      \end{multline}
			
			Explicitly, given a triple of tangent vectors
			\[
			    \bu=(u_1,u_2), \bv=(v_1,v_2), \bw=(w_1,w_2)\in  T_b\scB_{\fs\fo_5}=H^0(X,K_X^2)\oplus H^0(X,K_X^4)
			\]
			one has
			
			\begin{multline}
				c_b(\bu,\bv,\bw)=
				12\sum_{z_0\in V(b_2)} \frac{u_2 v_2 w_2}{b_1 j^1(b_2)^2}(z_0)+ \\
				24\sum_{z_0\in V(b_1^2-4b_2)}\frac{(b_1 u_1-2u_2)(b_1 v_1-2v_2)(b_1 w_1-2w_2)}{b_1 j^1(b_1^2-4b_2)^2}(z_0).
			\end{multline}

		\end{thm}

		We observe  that the formule (\ref{cubic_A2_fancy}), (\ref{cubic_B2_fancy}) and (\ref{cubic_G2_fancy})
		for the $A_2$, $B_2$ and  $G_2$ cubics have  similar  structure.
		In the $B_2$ and $G_2$ cases, when there are roots of different length, there are
		 two groups of summands, corresponding to the factorisation $\fD=\fD_1 \fD_2$.
		We expect that a similar formula holds in general and plan to discuss  this in a future work.

      \subsection*{Acknowledgements} P.D.~ thanks SISSA  for  hospitality.
	    
      \section{A brief reminder on $G_2$}\label{section_reminder_G2}

	  We begin with recalling some relevant facts about the $\fg_2$ root system,
	  and refer to   \cite{humphreys} or \cite{serre_Lie} for  systematic and beautiful expositions.

	The $\fg_2$ root system can be realized concretely as a subset $\Phi$ of eucliedan  $\RR^2$, 
	if one sets 
	 $\alpha_1=(1,0)$ and $\alpha_2=\left(-\frac{3}{2},\frac{\sqrt{3}}{2}\right)$. Then 
	 $\Phi\subseteq \RR^2$ consists of 
	 six short 
	roots, $\pm \alpha_1$, $\pm (\alpha_1+\alpha_2)$, $\pm (2\alpha_1+\alpha_2)$,
	having length $1$ and
	 six long ones, $\pm \alpha_2$, $\pm (3\alpha_1+\alpha_2)$, $\pm (3\alpha_1+2\alpha_2)$,  having length $\sqrt{3}$.

		      The root system is isomorphic to its dual root system, consisting of the coroots $\alpha^\vee=\frac{2\alpha}{\langle \alpha,\alpha\rangle}$, 
	$\alpha\in\Phi$.
		      An explicit isomorphism between $\Phi$ and $\Phi^\vee$ in the above realization is provided by the composition 
		      $\left(\textrm{mult}\frac{2}{\sqrt{3}}\right)\circ \left(\textrm{Rotation by }\frac{\pi}{2}\right): \Phi\to \Phi^\vee$.
		     We shall write  $\tau$ for the (anti-clockwise) rotation by $\frac{\pi}{2}$ in $\RR^2$  and notice that
		      \begin{equation}
			\tau: \alpha_1\longmapsto \frac{1}{\sqrt{3}}(3\alpha_1+ 2\alpha_2), \ \alpha_2\longmapsto -\sqrt{3}(2\alpha_1+\alpha_2).
		      \end{equation}

		     Instead of the euclidean metric we could have taken some scalar multiple of it (up to scaling, 
		     this is the only Weyl-invariant inner product), and that would affect the isomorphism between $\Phi$ and
		     $\Phi^\vee$.
		     In general, for two isomorphic root systems there is a linear isomorphism between the euclidean spaces
		     (containing them) that sends roots to roots and intertwines the Weyl group actions -- and such an
		     isomorphism is unique up to the action of the Weyl group(s).

		       We shall be working with the complexified root system, and treat it as the root system of the exceptional Lie group $G_2$, 
		       that is, $\Phi\subseteq \ft^\vee\simeq \CC^2$.

		     Let us recall how this is done concretely.

      The $14$-dimensional exceptional (complex) Lie group $G_2$ can be defined as the connected component of the automorphism group of
      the complex octonions (the Cayley algebra over $\CC$), see \cite{Schwarz_G2_spin7}. As such, it has a natural $7$-dimensional representation,
      the subspace of purely imaginary octonions. This representation is in fact orthogonal, preserving a naturally defined
      symmetric bilinear form. The Lie algebra $\fg_2$ is then identified with the Lie algebra of derivations of the octonions.
	  
	  It is possible to understand quite explicitly the $7$-dimensional representation  $\fg_2\subseteq \fs\fo_7$, see e.g. \cite[p.103]{humphreys},
	  \cite{adlervm_cis}.
	  In particular, the Cartan subalgebra   	 $\ft\subseteq \fg_2\subseteq \fs\fo_7$ consists
	  of diagonal matrices of the form
	    \[
		\ft= \left\{ h=\textrm{diag}(-a-b,-a,-b,0,b,a,a+b),\ a,b\in \CC \right\}.
	    \]
	   We can then take as simple roots the linear functionals $\alpha_1$, $\alpha_2\in \ft^\vee$ 
	\[
		\alpha_1\left(h\right)=  h_{55},\ \alpha_2\left(h\right)= h_{66}-h_{55}.
	\]
	
	The matrices in $\fg_2$   have  eigenvalues  $0, \pm\lambda_1,\pm \lambda_2,\pm \lambda_3$,
	$\sum_{i=1}^3 \lambda_i=0$, and one 
	can take as basic invariant polynomials $I_1$ and $I_2$ respectively
	$\frac{1}{2}(\lambda_1^2+\lambda_2^2+\lambda_3^2)$ and
	$(\lambda_1\lambda_2\lambda_3)^2$. This is discussed in  \cite{Schwarz_G2_spin7}, see also
	\cite{katz_pan_G2}, \cite{hitchin_G2}. Or,
	 one can 
	write out directly the characteristic polynomial of a matrix $h\in \ft$ as above, and obtain
	\[
		\det(h-\lambda E_7)= -\lambda^7+ \lambda^5 2I_1(h) -\lambda^3 I_1^2(h)  + \lambda I_2(h).
	\]
	The homogeneous invariant polynomials  generating $\CC[\ft]^W$ are not uniquely determined, but their degrees,
	$d_1=2$, $d_2=6$,
	are. In particular,  in some
	sources $I_1$ is defined without the $\frac{1}{2}$ factor.

	   For the calculation in   the paper we shall  fix a choice of simple roots
		       $\alpha_1,\alpha_2$ and use it  as a basis of $\ft^\vee$. Correspondingly, we  identify $\ft= \CC^2$ via the
		       (dual) basis.

		       Using the identifications $\ft=\CC^2$ and $\CC[\ft]^W = \CC[I_1,I_2]$, we can now identify the
		       adjoint quotient map $\chi: \ft\to \ft/W$ with the  $12:1$ map 
		       \begin{equation}
	  		\bI=(I_1,I_2): \CC^2\to \CC^2
			\end{equation}
		  where
		  \[
		  	\left|
			    \begin{array}{l}
			     I_1(\alpha_1,\alpha_2)= 3\alpha_1^2+3\alpha_1\alpha_2+\alpha_2^2\\
		  	I_2(\alpha_1,\alpha_2)= 4\alpha_1^6+12\alpha_1^5\alpha_2 + 13\alpha_1^4\alpha^2_2+ 6 \alpha_1^3\alpha_2^3+ \alpha_1^2\alpha_2^4.\\
			    \end{array}
		  	\right. 
		  \]

	   Next, the Weyl group of $G_2$ is is the dihedral group of order $12$, and since
		      $\tau\notin W=D_6$, while  $\tau^2=-1\in W$, we have that $\tau$ induces a \emph{non-trivial} involution $\tau$ of $\ft/W$.
	  We check with a direct calculation that 
			$\tau^\ast I_1=I_1$ and $\tau^\ast I_2= -I_2+ \frac{4}{27}I_1^3$.

	  We continue the abuse of notation and write $\tau: \cB\to \cB$ for the induced involution of the Hitchin base,
	   \begin{equation}\label{tau_base}
      	\tau (b_1,b_2)=\left(b_1,-b_2+\frac{4}{27}b_1^3\right).
		  \end{equation}
      
      We also write $\tau_i=\pr_i\circ \tau:\scB\to H^0(X,K_X^{d_i})$, $i=1,2$ for the two components of $\tau$,
      \begin{equation}\label{tau_components}
      	\tau_1=\pr_1,\ \tau_2=-\pr_2+\frac{4}{27}\pr_1^3.
      \end{equation}
      
        Our formulae for $\tau$ --  as an automorphism of $\ft$, $\ft/W$ or  $\cB$ -- depend on various choices.
	      For an intrinsic description, see Remark 1.3 in \cite{don-pan}. 
      
	      Up to scaling, the only Weyl-invariant symmetric bilinear form on $\ft$ is the Killing form $\kappa$ and
	      we know
	      \begin{equation}\label{killing}
	      	\kappa =\sum_{\alpha\in \Phi}\alpha\otimes \alpha\in \sym^2 \ft^\vee.
	      \end{equation}
	      With our choice of basis $\{\alpha_1,\alpha_2\}$, the Killing form is identified with the bilinear form
	      \[
	      	\kappa: \CC^2\times \CC^2\to \CC
	      \]
	      having   matrix $8\begin{pmatrix}
	      	                        	6&3\\
	      	                        	3&2\\
	      	                        \end{pmatrix}$
	      in the standard basis.

	      We also recall here that the discriminant of a Lie algebra  is the product of the roots,
	      $\fD=\prod_{\alpha\in \Phi}= (-1)^{|\Phi|/2}\prod_{\alpha\in \Phi^+}\alpha^2 $
	      which 
	       for $\fg_2$ is

	      \begin{equation}
	      \begin{split}
		\mathfrak{D}&= \alpha_1^2\alpha_2^2(\alpha_1+\alpha_2)^2(2\alpha_1+\alpha_2)^2(3\alpha_1+ \alpha_2)^2(3\alpha_1+2\alpha_2)^2\\
		&=I_2(4I_1^3-27 I_2)\in \CC[\ft]^W.
		\end{split}
	      \end{equation}
	      It is possible to obtain the last identty without too much brute force, for which we refer to  the discussion
	      surrounding  equation (36) in 
	      \cite{bruzzo_dalakov_sw}. In terms of the $\fs\fo_7$ matrices, this product is nothing but
	      $\lambda_1^2\lambda_2^2\lambda_3^2(\lambda_1-\lambda_2)^2(\lambda_1-\lambda_3)^2(\lambda_2-\lambda_3)^2$,
	      as, up to reordering, we can take $\alpha_1$ and $\alpha_2$ to be, respectively, $\lambda_1$ and $\lambda_2-\lambda_1$ (as above).
	      The two factors, $\fD_1=I_2$ and $\fD_2= 4I_1^3-27 I_2$, equal the  product of short and long roots, respectively.
	      (In general, we know that the Weyl group of a root system acts transitively on the set of roots of the same length.)
	      The second factor, $\fD_2$, is also nothing but $27\tau_2$, so
	       $\fD$ is not only Weyl-invariant, but it is also $\tau$-invariant.
	       
	       We are also going to write $\fD$ for the associated section of $\cO_{\cB}\otimes H^0(X,K_X^{12})$ -- or, if one prefers, a
	       global section
	       of the pullback of $K_X^{12}$ to $\cB\times X$. That is,

	  \begin{equation}\label{discriminant}
      	\fD=\fD_1\fD_2= 27 \pr_2\tau_2 \in \Gamma\left(\cB, \cO_\cB\ctimes H^0(X,K_X^{12})\right),
      	\end{equation}
      	\[ \fD(b_1,b_2)=b_2(-27b_2+4b_1^3)\in H^0(X,K_X^{12}).
	  \]

	\section{The Cameral Curve}\label{section_cam_curve}

	In this section we describe briefly the $G_2$ cameral curve and collect some of its features for future reference.
	
	Intrinsically,
	the $K_X$-valued cameral curve $\widetilde{X}_b$, corresponding to a  section $b\in \cB$
	is the pullback of the $W$-cover $\chi: \ft\ctimes K_X\to \ft\ctimes K_X/W$ by the evaluation map of the section
	$b\in H^0(X, \ft\ctimes K_X/W)$.
	With the choices discussed in section \ref{section_reminder_G2}, we have $\ft\ctimes K_X=K_X^{\oplus2}$, 
	$\ft\ctimes K_X/W=K_X^2\oplus K_X^6$ and
	$b=(b_1,b_2)\in \cB= H^0(X,K_X^2)\oplus H^0(X,K_X^6)$.
	We set $M=\textrm{tot}(K_X\oplus K_X)$ and write $\pi:M\to X$ for the bundle projection.

	Any linear map (e.g., any root) $\alpha\in\ft^\vee=\mhom(\ft,\CC)$ determines  a vector bundle homomorphism $\ft\ctimes K_X\to K_X$,
	denoted by the same letter -- and
	hence  a tautological section $\bal\in H^0(M, \pi^\ast K_X)$, mapping a (closed) point  $m\in M$ to 
	   $ \bal(m)=(m, \alpha(m))\in M\times_X \tot K_X$.

The cameral curve $\widetilde{X}_b\subset M$
	is cut out  by the equations

	\begin{equation} \label{G2_cameral}
	 	\left|
		  \begin{array}{l}
		  	3\bal_1^2+3\bal_1\bal_2+\bal_2^2=\pi^\ast b_1\\
		  	4\bal_1^6+12\bal_1^5\bal_2 + 13\bal_1^4\bal^2_2+ 6 \bal_1^3\bal_2^3+ \bal_1^2\bal_2^4= \pi^\ast b_2		\\
		  \end{array}
	 	\right. 
	 \end{equation}

	 which are sometimes more convenient in the form
	 \begin{equation} \label{G2_cameral_bis}
	 	\left|
		  \begin{array}{l}
		  	\frac{1}{2}\left\{ \bal_1^2+ (\bal_1+\bal_2)^2+(2\bal_1+ \bal_2)^2  \right\}  =\pi^\ast b_1\\
		  	\bal_1^2(\bal_1+\bal_2)^2(2\bal_1+\bal_2)^2 = \pi^\ast b_2.\\
		  \end{array}
	 	\right. 
	 \end{equation}

	 We  set $p_b=\left. \pi\right|_{\widetilde{X}_b}:\widetilde{X}_b\to X$.

	 The assignment $\alpha_i\mapsto \bal_i$ ($i=1,2$) determines an (injective) $\CC$-algebra homomorphism
		\begin{equation}
			\io: \sym(\bt^\vee)\hookr H^0\left(M, \bigoplus_{n\geq 0}\pi^\ast K_X^n\right),
		\end{equation}
	which we may occasionally suppress as well.

	  From the defining equations (\ref{G2_cameral}) or (\ref{G2_cameral_bis}) we have
	  
	  \begin{equation}\label{pullback_b2_first}
	  	 p_b^\ast(b_1)=\left. \io I_1\right|_{\widetilde{X}_b},\  p_b^\ast (b_2)= \left. \io I_2\right|_{\widetilde{X}_b},\  p_b^\ast(4b_1^3-27b_2)= \left. \io \left(4I_1^3-27I_2\right)\right|_{\widetilde{X}_b}.
	  \end{equation}

	The Balduzzi--Pantev formula (\ref{BP_ver1}), see also section \ref{section_periods}, involves certain quadratic residues at ramification points
	of $p_b:\widetilde{X}_b\to X$, so we describe the ramification points  explicitly.

	For $\alpha\in \ft^\vee$, the section $\bal$  vanishes along a 
	surface in the 3-fold $M$, the
	``hyperplane divisor'' 
		      $\tot \left(\ker \alpha\ctimes K_X\right) \subseteq M$,
		  a line subbundle of $\ft\ctimes K_X$.
		  The ramification $\textrm{Ram}(p_b)$ is determined by the intersection of $\widetilde{X}_b$ with the different root hyperplane divisors.
		Restricting  $\bal$  to $\widetilde{X}_b\subset M$ gives  a section
	$\bal_{\widetilde{X}_b}\in H^0(\widetilde{X}_b,p_b^\ast K_X)$.   We shall occasionally   suppress the subscript
	 indicating  the restriction. 
	    	The ramification divisor of $p_b$ is the sum of the different $D_{\alpha}=V(\bal_{\widetilde{X}_b})$, for
	 for $\alpha\in \Phi^+$, i.e., zeros of the sections arising from positive roots $\alpha$,
	 see \cite{scognamillo_elem}, \cite{don-gaits}. Notice it is enough to  restrict to positive roots only, as roots come in pairs.

	  The \emph{branch} divisor of $p_b$, on the other hand, is determined by the vanishing of $\fD(b)=b_2(4b_1^3-27b_2)$.
	  For a cameral curve with generic ramification, i.e., $b\in \scB$,
	  the zeros of $b_2\in H^0(X,K^6_X)$
	  and  $-27b_2+4b_1^3\in H^0(X,K_X^6)$ are simple and disjoint. 
	  The fibre over each branch point in this case contains  $|W|/2=6$ distinct points.

	  Using that the  positive roots of the $G_2$ root system are
	  $\alpha_1$,  $\alpha_1+\alpha_2$, $2\alpha_1+\alpha_2$ (short roots) and $3\alpha_1+\alpha_2$, $3\alpha_1+2\alpha_2$ and $\alpha_2$ (long roots), we
	  obtain 6 types of ramification points. All of them are cut out in $M$ by  equations  of the form
	  \[
	      \left|
		\begin{array}{l}
		L(\bal_1,\bal_2)=0\\
	  	I_1(\bal_1,\bal_2)=\pi^*b_1\\
	  	I_2(\bal_1,\bal_2)=\pi^*b_2\\
	  	\end{array}
	      \right. 
	  \]
	  where $L$ is the equation of the corresponding root hyperplane. Simplifying the respective systems a bit, we get:

		  \begin{align}\label{ram123}
			     (I) \begin{cases}\bal_1 = 0\\  \bal_2^2=\pi^*b_1\\ b_2=0\\  \end{cases} &
			     (II) \begin{cases}  \bal_1+\bal_2 = 0\\  \bal_1^2=\pi^*b_1\\  b_2=0\\ \end{cases} &
			     (III) \begin{cases}  2\bal_1+\bal_2 = 0\\  \bal_1^2=\pi^*b_1\\  b_2=0\\ \end{cases}
		  \end{align}
    and
		  
		  \begin{align}\label{ram456}
			     (IV) \begin{cases} 3\bal_1+\bal_2 = 0\\
			    \frac{1}{3}\bal_2^2=\pi^\ast b_1\\
			    \frac{4}{729}\bal_2^6=\pi^\ast b_2\\  \end{cases} &
			     (V) \begin{cases}   3\bal_1+2\bal_2 = 0\\
				    \frac{1}{3}\bal_2^2=\pi^\ast b_1\\
				    \frac{4}{729}\bal_2^6=\pi^\ast b_2\\ \end{cases} &
			     (VI) \begin{cases}  \bal_2 = 0\\
			    3\bal_1^2=\pi^*b_1\\
			    4\bal_1^6=\pi^*b_2\\   \end{cases}
		  \end{align}

	  Hence indeed  the short roots vanish  over (the support of)  $V(b_2)$  on $X$, while zeros of long
	  roots are above (the support of)  $V(\tau_2(b))=V(-27b_2+4b_1^3)$. Notice also that (by genericity) there are two 
	  distinct points, satisfying  each
	  of the above systems, so above each branch point there are indeed $2\times 3=6$ distinct ramification points.

		We also list here   the values of the Jacobi matrix $\io D\bI$ at the respective ramification points
		$x\in  D_\alpha$,
		
		\begin{align}\label{DI_ram}
		\io D\bI(x)=
		      \begin{cases}
		      	\begin{bmatrix}
	      	  3\bal_2& 2\bal_2\\
	      	  0& 0\\
	      	  \end{bmatrix}, & \alpha=\alpha_1 \\
	      	  \\
	      	  \begin{bmatrix}
				3\bal_1 & \bal_1\\
				0&0\\
			\end{bmatrix},  &  \alpha= \alpha_1+\alpha_2    \\
			\\
			\begin{bmatrix}
				0&-\bal_1\\
				0&0\\
			\end{bmatrix}, & \alpha= 2\alpha_1+\alpha_2  \\
			\\
			\begin{bmatrix}
				\bal_2 & \bal_2\\
				\frac{4}{81}\bal_2^5& \frac{4}{81}\bal_2^5\\
			\end{bmatrix}, & \alpha= 3\alpha_1+\alpha_2 \\
			\\
			\begin{bmatrix}
				-\bal_2& 0\\
				-\frac{4}{81}\bal_2^5&0\\
			\end{bmatrix}, & \alpha= 3\alpha_1+2\alpha_2 \\
			\\
			\begin{bmatrix}
				6\bal_1&3\bal_1\\
				24\bal_1^5& 12\bal_1^5\\
			\end{bmatrix}, & \alpha= \alpha_2.\\
		      \end{cases}
		\end{align}

		 \section{Infinitesimal Period Maps}\label{section_periods}

	      The first-order deformations of a polarized   abelian variety $Y$ are given by 
	      \[
	      	H^1(Y,T_Y)=H^1(Y, T_{Y,0}\ctimes \cO_Y)=T_{Y,0}\otimes H^1(Y,\cO_Y)= \otimes^2 T_{Y,0},
	      \]
	      where $T_{Y,0}$ is the tangent space to $Y$ at the origin. By Riemann's bilinear relations, the  deformations of $Y$ as an
	      abelian variety are in fact contained in the second symmetric power, $\textrm{Sym}^2T_{Y,0}$, of the tangent space of $Y$.
	      
	      If we have, more generally, a family of polarized abelian varieties $f:\cY\to \cB$, and $\cP$ is the period map of the associated
	      variation of integral weight-one polarized Hodge structures (\cite{griffiths_periods},\cite{voisin1}), then its differential is a bundle homomorphism
	      \begin{equation}
	      	 d\cP: T_{\cB}\to \sym^2 \cV,
	      \end{equation}
	      where $\cV=f_\ast (\ker df)$ is the vertical bundle of $f$.
	      While Hodge theory provides transcendental invariants of algebraic varieties, infinitesimal Hodge theory gives algebraic invariants. In particular, 
	      in many cases derivatives of period maps can be computed by residue formulae.
	      
	      Let us assume now that $\dim \cY=2\dim \cB$.
	      Unlike in the smooth setting, the total space $\cY$ of the  family of polarized abelian varieties need not carry a holomorphic symplectic structure for which the fibres
	      of $f$ are Lagrangian.
	      If such a Lagrangian structure does exist, however,  the symplectic form induces an isomorphism $i:\cV^\vee\simeq T_\cB$
	      in the familiar fashion. Donagi and Markman (\cite{donagi_markman_cubic}, \cite[\S 7.2]{donagi_markman}) show that if one fixes $f$ and an isomorphism $i$
	      as above, then there exists a Lagrangian structure on $\cY$ that induces $i$ if and only if $d\cP\circ i\in \mhom(\cV^\vee,\sym^2\cV)$
	      is actually in the image of $H^0(\cB, \sym^3\cV)\subseteq H^0(\cB, \cV^\vee\otimes \sym^2\cV)$. In fact, in the latter case there is a unique such
	      symplectic form for which the zero section is Lagrangian. In \cite[\S 7.2] {donagi_markman} the Poisson version of these statements is discussed
	      as well, where $i$ is not an isomorphism, but an injection -- but we do not need this level of generality here.
	      In the context of families of (gauged) Calabi--Yau threefolds, the
	       corresponding cubic has been known as the \emph{Yukawa cubic} or the \emph{Bryant--Griffiths cubic}, see \cite{bryant_griffiths}.

	      We are going to consider the family of generalized Prym varieties $f:\Prym_{\cX/\scB}\to\scB$ associated with the $G_2$ Hitchin system.
	      For Hitchin systems the Donagi--Markman cubic is given by the Balduzzi--Pantev formula \cite{balduzzi}, see also \cite{ugo_peter_cubic} for
	      the case of  generalized Hitchin systems. 
	      That formula states that, for $b\in\scB$,  $(d\cP\circ i)_b$ is the map 
	      \[
	      	H^0(\widetilde{X}_b,\ft\ctimes K_{\widetilde{X}_b})^W\to \sym^2\left(H^0(\widetilde{X}_b,\ft\ctimes K_{\widetilde{X}_b})^W\right)^\vee
	      \]
	      \begin{equation}
	      	\xi \longmapsto \left(\eta\cdot \zeta \mapsto  \frac{1}{2}\sum_{x\in \textrm{Ram}p_b}\textrm{Res}^2_x \left(p_b^\ast
	      	\frac{\cL_{i^{-1}(\xi)}\fD(b)}{\fD(b)} \eta\cup\zeta\right)\right).
	      \end{equation}
	      
	      Here $\cup$ is the cup product on sections of $\ft\ctimes K_{\widetilde{X}_b}$, which is obtained by combining the usual cup product on 
	      $K_{\widetilde{X}_b}$
	      with the Killing form on $\ft$. The first ingredient in deriving this  formula is   a result of Griffiths \cite[Theorem 1.27]{griffiths_periods},
	      relating
	       infinitesimal
	      period map to the cup product with the Kodaira--Spencer class of the family (here: of cameral Pryms).
	       The latter can be computed using  
	      the corresponding family of cameral curves $\cX\to \scB$. Their  deformations are essentially determined by varying the branch locus, enconded
	      in the discriminant.

	      The weight-one $\ZZ$-VHS arising from the cameral Pryms from above can be in fact identified with a weight-one $\ZZ$-VHS arising from $p:\cX\to \scB$,
	      see e.g. \cite{hertling_hoev_posthum}, \cite{don-pan}, \cite{Florian_HitchinCY}, 
	      \cite[\S 3.2]{bruzzo_dalakov_sw}. The latter has as its underlying flag
	      \[
		    R^0p_\ast^W\left(\ft\ctimes \Omega^1_{\cX/\scB}\right)\subseteq R^1 p_\ast^W\left(\ft\ctimes \Omega^\bullet_{\cX/\scB}\right),
	      \]
	      that is, to
	      $b\in \scB$ we associate the flag $H^0(\widetilde{X}_b,\ft\ctimes K_{\widetilde{X}_b})^W\subseteq H^1(\widetilde{X}_b, \ft)^W$.
	      The corresponding Gauss--Manin connection $\nabla^{GM}$ can be defined topologically or holomorphically, see 
	       \cite[\S 8] {hertling_hoev_posthum} for a discussion in the style of
	      \cite{katz-oda}. This variation of Hodge structures admits a Seiberg--Witten differential $\bla_{SW}$, see \cite[\S 8]{hertling_hoev_posthum},
	      \cite{Florian_HitchinCY}. The latter is a global section of $R^0p_\ast^W\left(\ft\ctimes \Omega^1_{\cX/\scB}\right)$, with value
	      $\bla_{SW}(b)= \left. \lambda\right|_{\widetilde{X}_b}\in H^0(\widetilde{X}_b,\ft\ctimes K_{\widetilde{X}_b})^W$, the restriction of
	      the $\ft$-valued Liouville form on $\textrm{tot}(\ft\ctimes K_{X})$ to $\widetilde{X}_b$.
	      
	      In this note we are going to use the isomorphism $\nabla^{GM}\bla_{SW}=i^{-1}:T_\scB\to \cV^\vee$, 
	      for which we have obtained an explicit description in \cite{bruzzo_dalakov_sw}. In terms of that, we have that the Balduzzi--Pantev formula for
	      \begin{equation}
	      	 c:= \sym^2(i^{-1})\circ d\cP\in H^0(\scB, \sym^3 T^\vee_\scB)\subseteq H^0(\scB, T^\vee\scB\otimes \sym^2 T^\vee_\scB)
	      \end{equation}
	      is given by the formula (\ref{BP_ver1}). Our goal will be to evaluate explicitly the quadratic residues appearing in this formula.

		\section{Residues And Contributions From The Discriminant}
		
		We begin evaluating the different contributions to the Balduzzi--Pantev formula in the form
		 (\ref{BP_ver1}).

	   First we compute the logarithmic derivative of the discriminant,
	 $\frac{\cL_u\fD}{\fD}$, as that enters the
      formula for the cubic. 
      From equation (\ref{discriminant}), we obtain
      
      \begin{equation}\label{log_der_discr}
      	\frac{d \fD}{\fD}=  \frac{d\fD_1}{\fD_1} + \frac{d\fD_2}{\fD_2}=\frac{d\pr_2}{\pr_2}  +  \frac{\frac{4}{9}\pr_1^2 d\pr_1 - d\pr_2}{\tau_2}. 
      \end{equation}
    At $b\in \cB$  the two summands in (\ref{log_der_discr}) have poles at $V(b_2)$ and $V(\tau_2(b))$, respectively.
      Explicitly, for  $b=(b_1,b_2)\in\scB$ and   $\bu=(u_1,u_2)\in T_b \scB=\cB$, this gives
      
      \begin{equation}\label{log_der_discr_eval}
      	\frac{d\fD_b(\bu)}{\fD(b)}= \frac{(d\pr_2)_b(\bu)}{\pr_2(b)}+\frac{(d\tau_2)_b(\bu)}{\tau_2(b)}
      	=\frac{u_2}{b_2}+\frac{-u_2+\frac{4}{9}b_1^2 u_1}{-b_2+\frac{4}{27}b_1^3}.
      \end{equation}
      
      We can now use the formulae (\ref{pullback_b2_first})  to explicate the pullback of this expression 
      by  $p_b$.

      \emph{Remark on notation: } 
       As
      $\cL_\bu \fD(b)= d\fD_b(\bu)$,  whether we  use $\cL$ or $d$ is mostly a matter of typography. 
      Both $\fD$ and $\cL_\bu\fD$ are holomorphic functions (on $\scB$) with values in a vector space (of dimension $>1$!), namely  $H^0(X,K_X^{12})$.
      Their ratio at $b\in \scB$ is understood as a meromorphic function on $X$, via  $\frac{\cL_\bu \fD}{\fD}(b)= \frac{\cL_\bu \fD(b)}{\fD(b)}$.
      Alternatively, one could think of $\frac{\cL_\bu\fD}{\fD}$ as a meromorphic function on $\scB\times X$ and restrict it to $\{b\}\times X$.

	 The formula (\ref{BP_ver1}) involves sums of terms which are 
	 quadratic residues of the meromorphic function $p_b^\ast\frac{\cL_\bu \fD}{\fD}(b)$, multiplied with an appropriate quadratic differential 
	 on $\widetilde{X}_b$. 
	 So, more generally, for a fixed $s\in H^0(\widetilde{X}_b,K_{\widetilde{X}_b}^2)$ we would like  to express in a more explicit way the residues
	 \begin{equation}
	 	\textrm{Res}^2_x\left(p_b^\ast \frac{d\fD_b(\bu)}{\fD(b)} s\right), 
	 \end{equation}
	 at  the various $x\in \textrm{Ram}(p_b)$. In order to shorten the argument, we are going to assume $s$ to be $W$-invariant, as it is 
	 in our intended applications, although
	 the result holds without that assumption as well.
	 
	 \emph{Remark:}\label{quad_res_def}
	We recall  that if $\zeta$ is a local coordinate on $\widetilde{X}_b$ centered at $x$, and if $p_b^\ast \frac{\cL_\bu \fD}{\fD}(b)s$ is
	 locally expressed as $\frac{\phi(\zeta)}{\zeta^2}d\zeta^2$, $\phi(0)\neq 0$, then 
	 $\textrm{Res}^2_x\left(p_b^\ast \frac{d_\bu \fD}{\fD}(b) s\right)=\phi(0)$.
	 
	 \begin{Prop}\label{residue_gen}
	 	After fixing  $b=(b_1,b_2)\in \scB$, let  $s\in H^0(\widetilde{X}_b, K^2_{\widetilde{X}_b})^W$  be a quadratic differential, and let 
	  $x\in D_\alpha= V(\bal_{\widetilde{X}_b})$ be a ramification point. Then
	
		\begin{multline} \label{residue_gen_formula}
\textrm{Res}^2_x\left(p_b^\ast \frac{d\fD_b(\bu)}{\fD(b)} \ s\right) = \\
			\begin{cases} \displaystyle 
				\textrm{Res}^2_x \left(\frac{1}{\bal^2} p_b^\ast\left(\frac{u_2}{b_1^2 }\right)s\right), 
				& \alpha \in \{\alpha_1,\alpha_1+\alpha_2, 2\alpha_1+\alpha_2\}\\[8pt]  \displaystyle
				\textrm{Res}^2_x\left(\frac{3}{\bal^2} p_b^\ast\left( \frac{ \frac{4}{9}b_1^2 u_1-u_2}{  b_1^2}\right) s\right),& 
				\alpha\in \{3\alpha_1+\alpha_2,3\alpha_1+2\alpha_2, \alpha_2\}\\
			\end{cases}
		\end{multline}
		that is, in terms of $\fD=\fD_1\fD_2$, 
		  
		  \begin{align} 
			\textrm{Res}^2_x\left(p_b^\ast \frac{d\fD_b(\bu)}{\fD(b)} \ s\right) =
			\begin{cases} \displaystyle
				\textrm{Res}^2_x\left(\frac{1}{\bal^2}	p_b^\ast \left(\frac{(d\fD_1)_b(\bu)}{b_1^2}\right) s\right),
				& \textrm{ for } \alpha \textrm{ short}\\[8pt] \displaystyle
				\frac{1}{9} \textrm{Res}^2_x\left(\frac{1}{\bal^2}	p_b^\ast \left(\frac{(d\fD_2)_b(\bu)}{ b_1^2}\right) s\right),
				 & \textrm{ for } 
				\alpha \textrm{ long}\\
			\end{cases}
		\end{align}

	  \end{Prop}
	  
	  Here and elsewhere we write $x\in D_\alpha$  instead of $x\in \textrm{supp }D_\alpha$, but this should not lead to confusion, 
	  as the ramification is generic, i.e., all $x\in \textrm{supp }D_\alpha$ have ramification index $2$.
	  
	  Notice that, by the genericity assumption, $p_b^\ast b_1(x)\neq 0$ -- so 
	  on the right hand side of  (\ref{residue_gen_formula})
	  the expression that is
	  pulled back from $X$ is a  quadratic differential, regular at $x$.

	    \emph{Proof:}

	 The statement can be verified case by case, or one can use the Weyl-invariance,  but the strategy is the same in all cases.  
	 The logarithmic derivative of $\fD$ splits as the sum of two terms, corresponding to products of short
	 and long roots, respectively, see (\ref{log_der_discr}). In each of these summands, 
	 we isolate in the denominator the  root that vanishes quadratically, and replace the
	 rest of the expression with the pullback of a section, having the same value at $x$ as the original expression, 
	 hence not affecting the quadratic residue.

	  We consider the short roots first. If $\bal$ is a short root, 
	 by the formula (\ref{log_der_discr})
	 we have 
	   \[
	  	 p_b^\ast \frac{d_b \fD(\bu)}{\fD(b)}=  p_b^\ast \frac{(d\fD_1)_b(\bu)}{\fD_1(b)} +\ldots   =p_b^\ast \frac{u_2}{b_2}+\ldots,
	  \]
	  a meromorphic function on $\widetilde{X}_b$, where the ellipsis represent
	   terms   regular at  $D_\alpha$,  so

	  \begin{equation}\label{res_simpl_1}
	 	\textrm{Res}^2_x\left(p_b^\ast \frac{d\fD_b(\bu)}{\fD(b)} s\right) = \textrm{Res}^2_x\left( \frac{p_b^\ast u_2}{p_b^\ast b_2} s\right)
	 	=
	 	\textrm{Res}^2_x\left( \frac{p_b^\ast u_2}{\bal_1^2(\bal_1+\bal_2)^2(2\bal_1+\bal_2)^2} s\right),
	 \end{equation}
	 using (\ref{pullback_b2_first}). Suppose $\bal=\bal_1$.
	 Then we have, using (\ref{ram123}),
	   \[
		  \frac{p_b^\ast b_2}{\bal_1^2}(x)= \left.(\bal_1+\bal_2)^2(2\bal_1+ \bal_2)^2\right|_{\bal_1=0}= \bal_2^4(x)=
		 p_b^\ast b_1^2(x)\neq 0,
	  \]
	  being   $b_1^2(p_b(x))$.
	  Hence we can rewrite the quadratic residue as 
	  
	  \begin{equation}
	  	\textrm{Res}^2_x\left(p_b^\ast \frac{d\fD_b(\bu)}{\fD(b)} s\right) =
	  	\textrm{Res}^2_x\left( \frac{1}{\bal_1^2} \frac{p_b^\ast u_2}{p_b^\ast b_1^2} s\right)=
	  	\textrm{Res}^2_x\left(\frac{1}{\bal^2}	p_b^\ast \left(\frac{(d\fD_1)_b(\bu)}{b_1^2}\right) s\right) .
	  \end{equation}

	  We can proceed  similarly with  the other short roots, using
	   (\ref{res_simpl_1}), (\ref{pullback_b2_first}) and (\ref{ram123}) again. 
	   However, this is not really needed, as the Weyl group acts transitively
	   on the set of roots of the same length, so the quadratic residues are the same.

	 For the ramification points corresponding to long roots, we proceed analogously.
	 This time we have that
	  by the formula (\ref{log_der_discr})
	  \begin{multline}
	  	p_b^\ast \frac{d\fD_b(\bu)}{\fD(b)} =  p_b^\ast \frac{(d\fD_2)_b(\bu)}{\fD_2(b)} +\ldots   =  p_b^\ast \frac{12b_1^2 u_1 - 27 u_2}{4b_1^3-27b_2} + \ldots 
	  	=  \\ \frac{p_b^\ast(12b_1^2 u_1 - 27 u_2)}{\bal_2^2 (3\bal_1+\bal_2)^2(3\bal_1+ 2\bal_2)^2}+ \ldots,
	  \end{multline}
	  and the terms in the elipsis do not contribute to the quadratic residue.
	  
	  If $\bal=3\bal_1+\bal_2$, we have
	  \[
	  	\frac{p_b^\ast(4b_1^3-27b_2)}{(3\bal_1+\bal_2)^2}(x)= \left. \bal_2^2(3\bal_1+2\bal_2)^2\right|_{3\bal_1+\bal_2=0}=\bal_2^4(x)=9p_b^*b_1(x),
	  \]
	  using (\ref{ram456}).
	  We can proceed similarly with the other long roots, or use the Weyl-invariance of the residues.

	  \qed
	  
	  In fact, the above proposition is also true without the requirement that $s$ be $W$-invariant, as one can see by going through 
	  all the choices for $\alpha$.

	  \section{Derivatives of $\bla_{SW}$ and Proof of Theorems \ref{thm1} and \ref{thm2} }\label{section_pf}

		  From  Proposition \ref{residue_gen} we see that in order to
		   evaluate the quadratic residues in (\ref{BP_ver1}), 
		   it is enough  to evaluate
	          $(\nabla^{GM}\bla_{SW})_b$ at ramification points. In our work \cite[Theorem A] {bruzzo_dalakov_sw}we have given
	          an expression for the latter, involving Lie-theoretic data, notably (the inverse of)  the Jacobi matrix
	          $D\bI$.

	           To make sense of evaluating the section $(\nabla^{GM}_\bv\bla_{SW})_b$ for $\bv\in T_b\scB$ at a  point we need a local trivialization
	           of $\ft\ctimes K_{\widetilde{X}_b}$ around that  point. For that
	           we are going to use local coordinates, as in \cite[\S 2.2]{bruzzo_dalakov_sw}.
	            We identify an open $U\subseteq X$
	           with a disk $\Delta\subseteq \CC$ and the bundle projection $\pi_U:M_U\to U$ with the canonical projection
	           $\Delta\times \CC^2\to \Delta$, with $\CC^2=\ft$ as before.
		    Away from the ramification locus we can use as a local coordinate
		    on $\widetilde{X}_U$ the pullback $p_b^\ast z$ of a local coordinate, $z$, on $U$. On the other hand, if we take $U$ to be centered at $z_0=p_b(x)\in \textrm{Bra}(p_b)$, we can use, as local coordinate on
		    $\widetilde{X}_U$ around $x$, one or both of $\alpha_1$ and $\alpha_2$. From the list in equation
		    (\ref{DI_ram}) one can see  which options are actually allowed by the implicit function theorem at the different   ramification
		    points.

		\begin{Prop}\label{nabla_lambda}
		
		Consider
	   	 $b=(b_1,b_2)\in \scB$,   $\bv=(v_1,v_2)\in T_b\scB=\cB_{\fg_2}$ and a ramification point 
	   	$x\in D_\alpha$ of $p_b:\widetilde{X}_b\to X$. Let
	   	$(U,\psi)$ be a local chart on $X$ and $z$ the respective
	   	   local coordinate,  centered at $z_0=p_b(x)$. Then

	   	\begin{align}
	   		(\psi^{-1})^\ast(\nabla^{GM}_\bv\bla_{SW})_b (x)=
			    \begin{cases}
			    	\frac{(d\fD_1)_b(\bv)[dz]}{j^1(\fD_1(b))}(z_0) [d\alpha]_x \ctimes
						    \begin{cases}
						           \frac{1}{2}
	   	 	      	\begin{pmatrix}
	   	 	      		2\\ -3\\
	   	 	      	\end{pmatrix} 
				      & \bal = \bal_1\\
				      \frac{1}{2}
	   	 	      	\begin{pmatrix}
	   	 	      		1\\ -3\\
	   	 	      	\end{pmatrix} 
				     & \bal= \bal_1+\bal_2\\
				     \frac{1}{2}
	   	 	      	\begin{pmatrix}
	   	 	      		1\\ 0\\
	   	 	      	\end{pmatrix}
				    & \bal= 2\bal_1+\bal_2\\ 
						    \end{cases} \\
						    \\
				\frac{(d\fD_2)_b(\bv)[dz]}{j^1(\fD_2(b))}(z_0) [d\alpha]_x \ctimes
					      \begin{cases}
						\frac{1}{2}
					         \begin{pmatrix}
	   	 	      		1\\ -1\\
	   	 	      	\end{pmatrix} & \bal =  3\bal_1+\bal_2\\
				      \frac{1}{2}
					   \begin{pmatrix}
	   	 	      		0\\ 1\\
	   	 	      	\end{pmatrix} & \bal=  3\bal_1+2\bal_2\\
				      \frac{1}{2}
				    \begin{pmatrix}
	   	 	      		1\\ -2\\
	   	 	      	\end{pmatrix}& \bal=\bal_2\\
					      \end{cases}
			    \end{cases}
	   	\end{align}
	    and, explicitly, 
		    
	   	\begin{equation}
		 \frac{(d\fD_1)_b(\bv)[dz]}{j^1(\fD_1(b))}(z_0)= \frac{v_2 [dz]}{j^1(b_2)}(z_0),
	\end{equation}
	\begin{equation}
		 \frac{(d\fD_2)_b(\bv)[dz]}{j^1(\fD_2(b))}(z_0)= \frac{(-v_2+\frac{4}{9}b_1^2v_1)[dz]}{j^1\left(-b_2+ \frac{4}{27}b_1^3\right)}(z_0).
	\end{equation}

	   	 \end{Prop}

	 We recall that in the first three cases $z_0\in V(b_2)$, while in the second three $z_0\in V(-27b_2+4 b_1^3)$,
	hence the 1-jets (at $z_0$) of $b_2$ and $-27b_2+4b_1^3$  can be considered as being in the fibre of $K_X^7$ (at $z_0$), see
	the equations (\ref{Atiyah_sqs}), (\ref{jet_canonical}).

	    For the expression of these quantities in local coordinates, see the proof.

	      \emph{Proof:}
	      
	      We proceed by applying  \cite[Theorem A]{bruzzo_dalakov_sw}, which states that
	      \begin{equation}
	      	(\nabla^{GM}_\bv\bla_{SW})_b = - \left. \io (D\bI)^{-1}\cdot \pi^\ast \bv\right|_{\widetilde{X}_b}.
	      \end{equation}
	      This is a global holomorphic section of $\ft\ctimes K_{\widetilde{X}_b}= K_{\widetilde{X}_b}^{\oplus 2}$, which  we 
	      will  evaluate  (using a trivialisation) at a ramification point
	      $x\in\textrm{Ram}(p_b)$, where the above formula has an apparent singularity.

	      Recall from \cite[\S 2]{bruzzo_dalakov_sw} that
	      the choice of local chart $(U,\psi)$ on $X$ and a compatible bundle chart $\phi: M_U\to \Delta\times \CC^2$
	      identifies $\widetilde{X}_U$ with a subset of $\Delta\times \CC^2$ and
	      the tangent bundle $T_{\widetilde{X}_U}$ with 
	      \[
		  \left. \ker \left(-\bbe' D\bI\right)\right|_{\phi\left(\widetilde{X}_U\right)}\subseteq 
		  \left. T_{\Delta\times\CC^2}\right|_{\phi\left(\widetilde{X}_U\right)},
	      \]
	      where $\bbe=\begin{pmatrix}
	                  	\beta_1\\ \beta_2\\
	                  \end{pmatrix}$ and $\beta_i:\Delta\to\CC$ are defined by $(\psi^{-1})^\ast b_i= \beta_i dz^{d_i}$, $i=1,2$.
	      Notice also that since we are using the chosen simple roots to identify $\ft=\CC^2$, the choice of bundle chart
	      $\phi$ is not an additional datum, but is determined by $\psi$.
	      
	      Let us denote  by $M_k$, $k=1,2$ the  $2\times 2$ matrix, obtained from $ \left(-\bbe' D\bI\right)$ by removing
	      the $(k+1)$-st column, that is, $\partial_k \bI$. The opens determined by $\det D\bI\neq 0$ and $\det M_k\neq 0$ cover
	      $\phi\left(\widetilde{X}_U\right)$, and on  these $[dz]$, respectively $[d\alpha_k]$ generates
	      the canonical bundle. 
	      We write $[dz]$ and $[d\alpha_k]$ instead of $dz$ and $d\alpha_k$ to emphasise that  the canonical bundle of 
	      $\widetilde{X}_U$ is a quotient of  $\left. \Omega^1_{M_U}\right|_{\widetilde{X}_U}$.
	      
	      A simple calculation shows, see \cite{bruzzo_dalakov_sw}, that on the intersection of two such opens
	      \begin{equation}
		    [dz]= \left. (-1)^{k+1}\frac{\det D\bI}{\det M_k}\right|_{\phi\left(\widetilde{X}_U\right)} [d\alpha_k].
	      \end{equation}
	      
	      Since locally $\bal_k$ is identified with $\alpha_k [dz]$, we obtain that $(\psi^{-1})^\ast(\nabla^{GM}_\bv\bla_{SW})_b$
	      equals, on $\{\det M_k\neq 0\}$, 
	      \begin{equation}\label{der_SW_expl}
	      	(-1)^k \left\{ \frac{\det[\bdelta\ \partial_2 \bI ]}{\det M_k}e_1 +  \frac{\det [\partial_1 \bI\ \bdelta]}{\det M_k} e_2  \right\} \otimes [d\alpha_k],
	      \end{equation}
	      where $\{e_1,e_2\}$ is the canonical basis on $\CC^2=\ft$ and $\bdelta= \begin{pmatrix}
	                                                                           	\delta_1\\ \delta_2\\
	                                                                           \end{pmatrix}$.
	      The functions $\delta_i:\Delta\to \CC$ are the coordinate expressions of the components of $\bv=(v_1,v_2)$, 
	      i.e., $(\psi^{-1})^\ast v_i= \delta_i dz^{d_i}$, $i=1,2$.

	      From  equation (\ref{DI_ram}), we see that for $x\in V(2\bal_1+\bal_2)$ only $\det M_1(\phi(x))\neq 0$, for 
	      $x\in V(3\bal_1+2\bal_2)$ only $\det M_2(\phi(x))\neq 0$, and in all other cases both $\det M_1$ and
	      $\det M_2$ are non-zero at $\phi(x)$.

	      We can compute all the needed determinants  using (\ref{DI_ram}) again. For example, 
	      
	      \begin{align}\label{detM}
		      \det M_1(\phi(x))=
		    \begin{cases}
		    	2\alpha_2 \beta_2'(\phi(x))=\frac{2}{3}\det M_2(\phi(x)),  & \alpha=\alpha_1 \\
		    	 (-3\alpha_1(4\alpha_1^4\beta_1'-\beta_2')(\phi(x))=\frac{1}{2}\det M_2(\phi(x)), & \alpha=\alpha_2 \\
		    \end{cases}
	      \end{align}
	      and
	      \begin{align}
		    \det [\bdelta\ \partial_2\bI](\phi(x)) =
		    \begin{cases}
		    	-2\alpha_2\delta_2(\phi(x)) =-\frac{2}{3}\det [\partial_1\bI\ \bdelta](\phi(x)), & \alpha=\alpha_1 \\
		    	 -3\alpha_1(4\alpha_1^4\delta_1-\delta_2)(\phi(x))= -\frac{1}{2}\det[\partial_1\bI\ \bdelta](\phi(x)),&\alpha=\alpha_2. \\
		    \end{cases}
	      \end{align}
	      
	      We can compute similarly the values of these expressions at the other roots -- or use the Weyl group action instead.
	      We then substitute  the last two equations into equation (\ref{der_SW_expl}). For instance, if 
	      $x\in D_{\alpha_1}$, 
	      \[
	      	(\psi^{-1})^\ast(\nabla^{GM}_\bv\bla_{SW})_b(x)= \frac{\delta_2}{\beta_2'}(z_0) \begin{pmatrix}
	      	                                                              	1\\ -3/2\\
	      	                                                              \end{pmatrix}\otimes [d\alpha_1]_x,
	      \]
	      and similarly for the other two short roots.
	      If
	      $x\in D_{\alpha_2}$,
	      \[
	      (\psi^{-1})^\ast(\nabla^{GM}_\bv\bla_{SW})_b(x)= -\frac{4\alpha_1^4\delta_1-\delta_2}{4\alpha_1^4\beta_1'-\beta_2'}(x)
	      	\begin{pmatrix}
	      		1\\ -2\\
	      	\end{pmatrix}\otimes [d\alpha_1]_x,
	      \]
	      and using (\ref{ram456}), we have that, as $\alpha_1^2(x)=\beta_1(z_0)/9$,
	      \[
	      	\frac{4\alpha_1^4\delta_1-\delta_2}{4\alpha_1^4\beta_1'-\beta_2'}(x)
	      	=\frac{-\delta_2+\frac{4}{9}\beta_1^2\delta_1}{-\beta'_2+\frac{4}{9}\beta_1^2\beta_1'}(z_0)
	      	=\frac{-\delta_2+\frac{4}{9}\beta_1^2\delta_1}{\left(-\beta_2+\frac{4}{27}\beta_1^3\right)'}(z_0).
	      \]
	      Similarly for the other two long roots.
	      
	      Notice that in the statement of the Proposition we have used $[d\alpha]_x$ as a generator of the fibre of $K_{\widetilde{X}_b}$
	      at $x$, and the previous formula is in terms of $[d\alpha_1]_x$.
	      We have from the earlier discussion 
	       that $[d\alpha_2]_x=-\frac{\det M_2}{\det M_1}(x) [d\alpha_1]_x$ on the respective opens.
	      For $\alpha=\alpha_2$ we have that $[d\alpha]_x= [d\alpha_2]_x=-2 [d\alpha_1]_x$, and this gives the stated result.
	      
	      Similarly for the other roots -- we can express $[d\alpha_1]_x$ or $[d\alpha_2]_x$ -- whichever we have used in
	      equation (\ref{der_SW_expl}) -- as a suitable multiple of $[d\alpha]_x$. Alternatively, one can use that
	      $(\nabla^{GM}_\bv \bla_{SW})_b$ is Weyl-invariant.

	      \qed

	      \emph{Proof of Theorem \ref{thm1}}

	      The formula for the Killing form (\ref{killing})  gives us that
	      \begin{equation}
	      	\kappa\left(
		      \begin{pmatrix}
		      	2\\ -3\\
		      \end{pmatrix},
		      \begin{pmatrix}
		      	2\\ -3\\
		      \end{pmatrix}
	      	\right)= 
	      	\kappa\left(
		      \begin{pmatrix}
		      	-1\\ 3\\
		      \end{pmatrix},
		      \begin{pmatrix}
		      	-1\\ 3\\
		      \end{pmatrix}
	      	\right)=
	      	 \kappa\left(
		      \begin{pmatrix}
		      	1\\ 0\\
		      \end{pmatrix},
		      \begin{pmatrix}
		      	1\\ 0\\
		      \end{pmatrix}
	      	\right)= 48,
	      \end{equation}
	      and
	      \begin{equation}
	      	\kappa\left(
		      \begin{pmatrix}
		      	-1\\ 1\\
		      \end{pmatrix},
		      \begin{pmatrix}
		      	-1\\ 1\\
		      \end{pmatrix}
	      	\right)= 
	      	\kappa\left(
		      \begin{pmatrix}
		      	0\\ 1\\
		      \end{pmatrix},
		      \begin{pmatrix}
		      	0\\ 1\\
		      \end{pmatrix}
	      	\right)=
	      	 \kappa\left(
		      \begin{pmatrix}
		      	1\\ -2\\
		      \end{pmatrix},
		      \begin{pmatrix}
		      	1\\ -2\\
		      \end{pmatrix}
	      	\right)= 16.
	      \end{equation}
	      
	      Consequently, from Proposition \ref{nabla_lambda} we obtain that at $x\in D_\alpha$
	      \begin{multline}
		      (\nabla^{GM}_\bv\bla_{SW})_b\cup (\nabla^{GM}_\bw\bla_{SW})_b (x)=\\
		      = \begin{cases} \displaystyle
				  12 \frac{(d\fD_1)_b(\bv)(d\fD_1)_b(\bw)[dz]^2}{j^1(\fD_1(b))^2}(z_0) [d\alpha]^2_x& \alpha \textrm{ short}\\[8pt] \displaystyle
				  4 \frac{(d\fD_2)_b(\bv)(d\fD_2)_b(\bw)[dz]^2}{j^1(\fD_2(b))^2}(z_0)[d\alpha]^2_x& \alpha \textrm{ long}\\
		        \end{cases}\\
			 =\begin{cases} \displaystyle
			  	12 \frac{v_2 w_2[dz]^2}{j^1(b_2)^2}(z_0)[d\alpha]^2_x & \alpha \textrm{ short }\\[8pt] \displaystyle
			  	4
				    \frac{(-v_2+\frac{4}{9}b_1^2v_1) (-w_2+\frac{4}{9}b_1^2w_1)[dz]^2}{j^1\left(-b_2+ \frac{4}{27}b_1^3\right)^2}(z_0)[d\alpha]^2_x  
				    & \alpha \textrm{ long } 
			  \end{cases}
	   	 \end{multline}

	      Now from  Proposition \ref{residue_gen}  it follows that the residues are
	      \begin{multline}\label{residues_final_global}
	      	\textrm{Res}^2_x\left(p_b^\ast \frac{d_b \fD(\bu)}{\fD(b)} (\nabla^{GM}_\bv\bla_{SW})_b\cup (\nabla^{GM}_\bw\bla_{SW})_b  \right) =\\
	      	=\begin{cases} \displaystyle
			   12 \frac{(d\fD_1)_b(\bu)(d\fD_1)_b(\bv)(d\fD_1)_b(\bw)}{b_1^2 j^1(\fD_1(b))^2}(z_0) & \alpha \textrm{ short}\\[8pt] \displaystyle
				  \frac{4}{9} \frac{(d\fD_2)_b(\bu)(d\fD_2)_b(\bv)(d\fD_2)_b(\bw)}{b_1^2 j^1(\fD_2(b))^2}(z_0)& \alpha \textrm{ long}\\ 
	      	 \end{cases}\\
	      	=
	      	\begin{cases} \displaystyle
	      		12\frac{u_2 v_2 w_2}{b_1^2 j^1 (b_2)^2}(z_0),& \alpha \textrm{ short}\\[8pt] \displaystyle
	      		12 \frac{(-u_2+\frac{4}{9}b_1^2 u_1)(-v_2+\frac{4}{9}b_1^2v_1) (-w_2+\frac{4}{9}b_1^2w_1)}{b_1^2j^1\left(-b_2+ \frac{4}{27}b_1^3\right)^2}(z_0)
	      		& \alpha \textrm{ long} \\
	      	\end{cases}.
	      \end{multline}

	      Now, for a branch point $z_0$, $\left| p^{-1}_b(z_0)\right|=6$.
	      The support of each divisor $D_\alpha= V(\bal_{\widetilde{X}_b})$ 
	      contains \emph{two} distinct ramification points (at which the residues are the same).
	      Residues are the same for roots of the same length, and there are three (positive) roots of each kind. 
	      But the formula (\ref{BP_ver1}) contains an overall factor of $\frac{1}{2}$, 
	      hence we obtain the overall factor of $36=12\cdot 3\cdot 2\cdot \frac{1}{2}$
	      from the statement of the 
	      Theorem.

		\qed

	      We mention for completeness that one can easily write the residues
	          in terms of a local coordinate: setting    $b_i=\beta_i(z)dz^{d_i}$,
	         $u_i=\gamma_i(z)dz^{d_i}$,
	         $v_i=\delta_i(z)dz^{d_i}$,
			$w_i=\epsilon_i(z) dz^{d_i}$,  $i=1,2$ we have
	      \begin{multline}\label{residues_final_local}
	      	\textrm{Res}^2_x\left(p_b^\ast \frac{d_\bu \fD}{\fD}(b) (\nabla^{GM}_\bv\bla_{SW})_b\cup (\nabla^{GM}_\bw\bla_{SW})_b  \right) =\\
	      	\begin{cases} \displaystyle
	      		12\frac{\gamma_2 \delta_2 \epsilon_2}{\beta_1^2  (\beta'_2)^2}(z_0),& \alpha \in \{\alpha_1,\alpha_1+\alpha_2, 2\alpha_1+\alpha_2\}  \\[8pt] \displaystyle
	      		12 \frac{(-\gamma_2+\frac{4}{9}\beta_1^2 \gamma_1)(-\delta_2+\frac{4}{9}\beta_1^2\delta_1) (-\epsilon_2+\frac{4}{9}\beta_1^2\epsilon_1)}{\beta_1^2\left(-\beta_2+ \frac{4}{27}\beta_1^3\right)^{'2}}(z_0),
	      		& \alpha\in \{3\alpha_1+\alpha_2,3\alpha_1+2\alpha_2, \alpha_2\}   \\
	      	\end{cases}.
	      \end{multline}

	      \emph{Proof of Theorem \ref{thm2}:}
	      
	      The expression for $c$ in the statement of Theorem \ref{thm2} is clear once we recall 
	      $\fD_2=27\tau_2$ and
	      that we have already
	      evaluated $d\tau_2$ in equation (\ref{log_der_discr}), see also (\ref{log_der_discr_eval}).
	      
	      This immediately implies the $\tau$-invariance of $c$, since the latter means  that $\tau^\ast c=c$, and we have that by definition
	      \begin{equation}
	      	(\tau^\ast c)_b(\bu,\bv,\bw)= c_{\tau(b)}(d\tau_b(\bu), d\tau_b(\bv), d\tau_b(\bw)).
	      \end{equation}
	      The involution $\tau$ exchanges the two summands in Theorem \ref{thm1}: $\pr_2(\tau(b))=\tau_2(b)$, $\pr_2\tau_2(\tau(b))=\pr_2(b)=b_2$
	      and ditto for the differentials $d\pr_2$ and $d\tau_2$.
	      
	      \qed

	      \section{Appendix: the $A_2$ case}\label{appendix_A2}
	      
	      Here we sketch the calculation of the Donagi--Markman cubic for the case of Hitchin systems of type $A_2=\fs\fl_3(\CC)$.
	      The Hitchin base
	      in this case is $\cB_{\fs\fl_3}=H^0(X, K_X^2)\oplus H^0(X,K_X^3)$, and we are going to write $\scB\subseteq \cB$ for the Zariski open set of
	      generic cameral covers. 
	      
	      We recall (e.g., \cite[\S 5.2]{bruzzo_dalakov_sw}) that we can identify the adjoint quotient  with the map
	      \begin{equation}
		    \bI:\CC^2\to \CC^2,
	      \end{equation}
	      
	      \begin{equation}\label{sl3}
	  	\left|
		    \begin{array}{l}
		    I_1(\alpha_1,\alpha_2)= \alpha_1^2 +\alpha_1\alpha_2 +\alpha_2^2\\
		    I_2(\alpha_1,\alpha_2)= -2\alpha_1^3-3\alpha_1^2\alpha_2 +3\alpha_1\alpha_2^2+ 2\alpha_2^3.\\
		    \end{array}
		\right. 
		\end{equation}
		Here the Cartan subalgebra $\ft\subseteq \fs\fl_3(\CC)$ consists of diagonal traceless $3\times 3$ matrices and we identify
		$\ft=\CC^2$ by choosing as simple roots $\alpha_1(A)=A_{11}-A_{22}$ and $\alpha_2(A)=A_{22}-A_{33}$. 
		There are three simple roots, all of the same length: $\alpha_1$, $\alpha_2$, $\alpha_1+\alpha_2$.
		The Weyl group is
		$W=S_3$, and we identify $\ft/W=\CC^2$ by the choice $I_1(A)= -3(A_{11}A_{22}+ A_{11}A_{33}+ A_{22}A_{33})$, 
		$I_2(A)=-27\det A$. 
		
		For  the $A_2$ root system, the discriminant is

	      \begin{equation}
	      \begin{split}
		\mathfrak{D}&= -\alpha_1^2\alpha_2^2 (\alpha_1+\alpha_2)^2\\ 
		&=-\frac{1}{27}(4I_1^3-I_2^2)\in \CC[\ft]^W.
		\end{split}
	      \end{equation}
	      
	      As before, we denote by $\fD$ also the associated section of $\cO_\cB\otimes H^0(X, K_X^6)$, i.e.
	      \begin{equation}\label{discriminant_A2}
      	\fD=-\frac{1}{27} \left(4\pr_1^3-\pr_2^2 \right)  \in \Gamma\left(\cB, \cO_\cB\ctimes H^0(X,K_X^{6})\right),
      	\end{equation}
      	\[ \fD(b_1,b_2)=-\frac{1}{27}(4b_1^3-b_2^2)\in H^0(X,K_X^{6}).
	  \]
	  
	  The logarithmic derivative of the discriminant is then seen to be
		\begin{equation}
		      \frac{d\fD}{\fD}= \frac{12\pr_1^2 d\pr_1-2\pr_2 d\pr_2}{4\pr_1^3-pr_2^2}.
		\end{equation}
		 That is, for $b=(b_1,b_2)\in\scB$ and  $\bu=(u_1,u_2)\in T_b \scB=\cB$,  
      
      \begin{equation}\label{log_der_discr_eval_A2}
	    \frac{\cL_\bu \fD}{\fD}(b)= \frac{12b_1^2 u_1-2b_2 u_2}{4b_1^3-b_2^2}.
      \end{equation}

	  The cameral cover $\widetilde{X}_b\subseteq \tot (K_X\oplus K_X)$ is cut out by
	  \begin{equation}\label{sl3_cameral}
		\left|
		    \begin{array}{l}
		    	\bal_1^2+\bal_1\bal_2+\bal_2^2=\pi^\ast b_1\\
		    	-2\bal_1^3-3\bal_1^2\bal_2+3\bal_1\bal_2^2+2\bal_2^3= \pi^\ast b_2\\
		    \end{array}
		\right. 
	\end{equation}
	and the   ramification divisors  $D_\alpha= V(\bal_{\widetilde{X}_b})$ of the 6-fold cover $p_b:\widetilde{X}_b\to X$
	  are determined by
	  \begin{align}\label{ram_A2}
			     (I) \begin{cases}\bal_1 = 0\\  \bal_2^2=\pi^*b_1\\ 2\bal_2^3=\pi^\ast b_2\\  \end{cases} &
			     (II) \begin{cases}  \bal_2 = 0\\  \bal_1^2=\pi^*b_1\\  -2\bal_1^3=\pi^\ast b_2\\ \end{cases} &
			     (III) \begin{cases}  \bal_1+\bal_2 = 0\\  \bal_1^2=\pi^*b_1\\  2\bal_1^3=\pi^\ast b_2\\ \end{cases}.
		  \end{align}
	  We  assume genericity, so the zeros of $b_1$ and $b_2$ are disjoint -- and hence disjoint from those of $\fD(b)$. 
	  
	  The values of the Jacobi matrix $\io D\bI$ at ramification points $x\in D_\alpha$ are now seen to be
	  \begin{align}\label{DI_ram_A2}
		\io D\bI(x)=
		      \begin{cases}
		      	\begin{bmatrix}
	      	  \bal_2& 2\bal_2\\
	      	  3\bal_2^2& 6\bal_2^2\\
	      	  \end{bmatrix}, & \alpha = \alpha_1 \\
	      	  \\
	      	  \begin{bmatrix}
				2\bal_1 & \bal_1\\
				-6\bal_1^2& -3\bal_1^2\\
			\end{bmatrix},  &  \alpha=\alpha_2   \\
			\\
			\begin{bmatrix}
				\bal_1&-\bal_1\\
				3\bal_1^2& -3\bal_1^2\\
			\end{bmatrix}, & \alpha=\alpha_1+ \alpha_2. \\
		      \end{cases}
		\end{align}
		
		There is a simplified formula for the residues in this case, as well: the analogue of Proposition \ref{residue_gen} is
		
		\begin{Prop}\label{residue_gen_A2}
			Let   $b=(b_1,b_2)\in \scB$, let  $s\in H^0(\widetilde{X}_b, K^2_{\widetilde{X}_b})^W$  be a quadratic differential,  and let 
	 	$x\in D_\alpha$ be a ramification point. Then

		\begin{multline}\label{residue_gen_formula_A2}
			\textrm{Res}^2_x\left(p_b^\ast \frac{d\fD_b(\bu)}{\fD(b)} \ s\right) =
			  -\frac{1}{27}\textrm{Res}^2_x \left(\frac{1}{\bal^2}p_b^\ast\left( \frac{12b_1^2 u_1-2b_2u_2}{b_1^2}\right) s \right)\\
			  =
			   \textrm{Res}^2_x\left(\frac{1}{\bal^2}p_b^\ast \left(\frac{d\fD_b(\bu)}{b_1^2}\right)s \right).
		\end{multline}
		\end{Prop}

		Once again, the same result is also true  without requiring the section $s$ to be Weyl-invariant.
		
		We next evaluate  $(\nabla^{GM}_\bv\bla_{SW})_b$ at a ramification point, using a local trivialization.

		\begin{Prop}\label{nabla_lambda_A2}
		Consider
	   	 $b=(b_1,b_2)\in \scB$,   $\bv=(v_1,v_2)\in T_b\scB=\cB_{\fs\fl_3}$ and a ramification point 
	   	$x\in D_\alpha$ of $p_b:\widetilde{X}_b\to X$. Let
	   	$(U,\psi)$ be a local chart on $X$ and $z$ the respective
	   	   local coordinate,  centered at $z_0=p_b(x)$. Then

	   	 \begin{align}
		      (\psi^{-1})^\ast(\nabla^{GM}_\bv\bla_{SW})_b (x)=
		     \frac{d\fD_b(\bv)[dz]}{j^1(\fD(b))}(z_0)\ [d\alpha]_x\ctimes
		      \begin{cases}
		      	\frac{1}{2}\begin{pmatrix}
		      	           	2\\ -1\\
		      	           \end{pmatrix}, & \alpha=\alpha_1 \\
		      	           \frac{1}{2}
			\begin{pmatrix}
				-1\\ 2\\
			\end{pmatrix},& \alpha=\alpha_2\\
			\frac{1}{2}
			\begin{pmatrix}
				1\\ 1\\
			\end{pmatrix}, & \alpha=\alpha_1+\alpha_2\\
		      \end{cases}
	   	 \end{align}
		  and, explicitly,
		  \begin{equation}
		  	\frac{d\fD_b(\bv)[dz]}{j^1(\fD(b))}(z_0)=  \frac{\left(12 b_1^2 v_1-2b_2 v_2 \right)[dz]}{j^1(4b_1^3-b_2^2)}(z_0).
		  \end{equation}

	   	 \end{Prop}

	   	\emph{Proof of Theorem \ref{thm3}:}
	   	
	   	 The Killing form $\kappa$ has matrix $2\begin{pmatrix}
	   	                                        	2& 1\\
	   	                                        	1& 2\\
	   	                                        \end{pmatrix}$ with respect to our chosen basis and
	   	    the vectors    $\frac{1}{2}\begin{pmatrix}
		      	           	2\\ -1\\
		      	           \end{pmatrix}$,
		      $ -\frac{1}{2}
			\begin{pmatrix}
				1\\ -2\\
			\end{pmatrix}$ and $\frac{1}{2}
			\begin{pmatrix}
				1\\ 1\\
			\end{pmatrix}$ all square to $3$ with respect to $\kappa$.

		We then use Proposition \ref{nabla_lambda_A2} to obtain that	
		\begin{multline}
	   	          (\psi^{-1})^\ast  (\nabla^{GM}_\bv\bla_{SW})_b\cup (\nabla^{GM}_\bw\bla_{SW})_b (x)=
	   	            3 \frac{d\fD_b(\bv) d\fD_b(\bw)[dz]^2}{j^1(\fD(b))^2}(z_0)\ [d\alpha]^2_x\\
		    =3  \frac{\left(12 b_1^2 v_1-2b_2 v_2 \right)\left(12b_1^2 w_1-2b_2 w_2 \right)[dz]^2}{j^1(4b_1^3-b_2^2)^2}(z_0)\ [d\alpha_1]^2_x                            
	   	 \end{multline}                          
	   	 
	   	 We then obtain, using Proposition \ref{residue_gen_A2}, that the residue at $x\in D_\alpha$ is

	   	 \begin{multline}
	   	 		\textrm{Res}^2_x\left(p_b^\ast \frac{d_b \fD(\bu)}{\fD(b)} (\nabla^{GM}_\bv\bla_{SW})_b\cup (\nabla^{GM}_\bw\bla_{SW})_b  \right) =
	   	 		3 \frac{d\fD_b(\bu) d\fD_b(\bv)d\fD_b(\bw)}{j^1(\fD(b))^2}(z_0)\\
	   	 		=-\frac{1}{9}\frac{\left(12 b_1^2 u_1-2b_2 u_2 \right)\left(12 b_1^2 v_1-2b_2 v_2 \right)\left(12b_1^2 w_1-2b_2 w_2 \right)}
	   	 		{j^1(4b_1^3-b_2^2)^2}(z_0).
	   	 \end{multline}
	   	 We then obtain the final formula from Theorem \ref{thm3} by taking into account that there are 3 simple roots and  $\left|p_b^{-1}(z_0)\right|=3$.

	   	 \qed

		\section{Appendix: the  $B_2$ case }\label{appendix_B2}
		
		We consider now the case of  $B_2=\fs\fo_5(\CC)$ Hitchin systems.
		
		The Hitchin base is $\cB_{\fs\fo_5}= H^0(X,K_X^2)\oplus H^0(X, K_X^4)$, and we write $\scB$ for
		the Zariski open subset  of generic $B_2$ cameral covers, reviewed below.
		
		We are going to consider $\fs\fo_5(\CC)$ as the Lie algebra
		of $5\times 5$ matrices, orthogonal with respect to the $5\times 5$ ``exchange matrix'' $J$ (having $1$'s along the main
		anti-diagonal and $0$'s elsewhere). Thus $A\in \fg\fl_5(\CC)$ is in $\fs\fo_5(\CC)$ when $A^TJ=-JA$. The Cartan subalgebra $\ft$ consists of
		diagonal matrices, satisfying $A_{11}=-A_{55}$, $A_{22}=-A_{44}$, $A_{33}=0$. A choice of simple roots 
		$\{\alpha_1, \alpha_2\}$ is  $\alpha_1(A)=A_{11}-A_{22}$, $\alpha_2(A)=A_{22}$. Then the corresponding positive roots are
		$\alpha_1$, $\alpha_1+\alpha_2$, $\alpha_1+2\alpha_2$, $\alpha_2$. Of these, $\alpha_1$ and $\alpha_1+2\alpha_2$ are long,
		while $\alpha_1+\alpha_2$ and $\alpha_2$ are short. 
		
		The Weyl group is $W= D_4$.
		A matrix in $B_2$ has eigenvalues $\pm \lambda_1$, $\pm \lambda_2$ and $0$, and the two Weyl-invariant quantities can be taken to be
		the sum and product of $\lambda_1^2$ and $\lambda_2^2$. 
		
		With these conventions, we identify the adjoint quotient map with
		
		\begin{equation}
			\bI: \CC^2\to \CC^2,
		\end{equation}
		
		\begin{equation}\label{so5}
	  	\left|
		    \begin{array}{l}
		    I_1(\alpha_1,\alpha_2)= \alpha_1^2 +2\alpha_1\alpha_2 +2\alpha_2^2\\
		    I_2(\alpha_1,\alpha_2)= \alpha_1^2\alpha_2^2+ 2\alpha_1\alpha_2^3+ \alpha_2^4,\\
		    \end{array}
		\right. 
		\end{equation}
		and these two polynomials can be also written as $(\alpha_1+\alpha_2)^2+\alpha_2^2$ and 
		$(\alpha_1+\alpha_2)^2\alpha_2^2$, respectively.

	        For the $B_2$ root system the discriminant is 

	      \begin{equation}
	      \begin{split}
		\mathfrak{D}&= \alpha_1^2(\alpha_1+2\alpha_2)^2\alpha_2^2 (\alpha_1+\alpha_2)^2\\ 
		&= I_2 (I_1^2-4I_2)\in \CC[\ft]^W.
		\end{split}
	      \end{equation}
	      
	      We write $\fD_1$ and $\fD_2$ for the two factors, $I_2$ and $I_1^2-4I_2$ respectively, corresponding to the products
	      of short and long roots, respectively.
	      
	      As before, we denote by $\fD$ also the associated section of $\cO_\cB\otimes H^0(X, K_X^4)$, i.e.
	      \begin{equation}\label{discriminant_B2}
      	\fD= \fD_1\fD_2= \pr_2(\pr_1^2-4\pr_2)  \in \Gamma\left(\cB, \cO_\cB\ctimes H^0(X,K_X^{4})\right),
      	\end{equation}
      	\[ \fD(b_1,b_2)= b_2(b_1^2-4b_2) \in H^0(X,K_X^{4}).
	  \]

		   The logarithmic derivative of the discriminant is then 
		\begin{equation}
		      \frac{d\fD}{\fD}= \frac{d\fD_1}{\fD_1} + \frac{d\fD_2}{\fD_2}= \frac{d\pr_2}{\pr_2} + \frac{2\pr_1 d\pr_1 - 4d\pr_2}{\pr_1^2-4\pr_2}.
		\end{equation}
		 That is, for $b=(b_1,b_2)\in\scB$ and  $\bu=(u_1,u_2)\in T_b \scB=\cB_{\fs\fo_5}$,  
      
      \begin{equation}\label{log_der_discr_eval_B2}
	    \frac{\cL_\bu \fD}{\fD}(b)= \frac{u_2}{b_2} + \frac{2b_1 u_1-4u_2}{b_1^2-4b_2}.
      \end{equation}

		   The cameral cover $\widetilde{X}_b\subseteq \tot (K_X\oplus K_X)$ is cut out by
	  \begin{equation}\label{so5_cameral}
		\left|
		    \begin{array}{l}
		    	\bal_1^2+2\bal_1\bal_2+2\bal_2^2=\pi^\ast b_1\\
		       \bal_1^2\bal_2^2+2\bal_1\bal_2^3+ \bal_2^4= \pi^\ast b_2\\
		    \end{array}
		\right. 
	\end{equation}
	and the   ramification divisors  $D_\alpha= V(\bal_{\widetilde{X}_b})$ of the cover $p_b:\widetilde{X}_b\to X$
	  are determined by
	  \begin{align}\label{ram_B2}
			      \begin{cases}\bal_1 = 0\\  2\bal_2^2=\pi^*b_1\\ \bal_2^4=\pi^\ast b_2\\  \end{cases} &
			      \begin{cases}\bal_1+2\bal_2=0\\ 2\bal_2^2=\pi^*b_1\\ \bal_2^4=\pi^\ast b_2\\   \end{cases}&
			      \begin{cases}  \bal_2 = 0\\  \bal_1^2=\pi^*b_1\\  \pi^\ast b_2=0\\ \end{cases} &
			      \begin{cases}  \bal_1+\bal_2 = 0\\  \bal_1^2=\pi^*b_1\\  \pi^\ast b_2=0\\  \end{cases}.
		  \end{align}
		  The support of each $D_\alpha$ contains two distinct points.
	  We  assume genericity, so the zeros of $b_1$ and $b_2$ are disjoint. The $8$-fold cover $p_b:\widetilde{X}_b\to X$ is branched over
	  $V(b_2)$ and $V(b_1^2-4b_2)$, and for each branch point $z_0$, $\left| p_b^{-1}(z_0)\right|=4$.
	  
	  The values of the Jacobi matrix $\io D\bI$ at ramification points $x\in D_\alpha$ are now seen to be, respectively
	  \begin{align}\label{DI_ram_B2}
		\io D\bI(x)=
		      \begin{cases}
		      	\begin{bmatrix}
	      	  2\bal_2& 4\bal_2\\
	      	  2\bal_2^3& 4\bal_2^3\\
	      	  \end{bmatrix}, & \alpha = \alpha_1 \\
	      	   & \\
	      	  \begin{bmatrix}
	      	  	-2\bal_2& 0\\
	      	  	-2\bal_2^3& 0\\
	      	  \end{bmatrix}, & \alpha=\alpha_1+2\alpha_2\\
	      	   & \\
	      	  \begin{bmatrix}
				2\bal_1 & 2\bal_1\\
				0 & 0 \\
			\end{bmatrix},  &  \alpha=\alpha_2   \\
			& \\
			\begin{bmatrix}
				0 & 2\bal_2\\
				0&0\\
			\end{bmatrix}, & \alpha=\alpha_1+ \alpha_2 \\
		      \end{cases}
		\end{align}

	We now state the analogue of Propositions \ref{residue_gen} and \ref{residue_gen_A2}	
	  \begin{Prop}\label{residue_gen_B2}
	 	Consider  $b=(b_1,b_2)\in \scB$, a quadratic differential $s\in H^0(\widetilde{X}_b, K^2_{\widetilde{X}_b})^W$ and
	 	a ramficiation point $x\in D_\alpha$. Then

		\begin{align}\label{residue_gen_formula_B2}
			\textrm{Res}^2_x\left(p_b^\ast \frac{d_b \fD(\bu)}{\fD(b)} \ s\right) =
			\begin{cases}
				\textrm{Res}^2_x \left(\frac{1}{\bal^2} p_b^\ast\left(\frac{u_2}{b_1 }\right)s\right),
				& \alpha\in \{\alpha_2, \alpha_1+\alpha_2\} \\
				\textrm{Res}^2_x\left(\frac{1}{\bal^2} p_b^\ast\left( \frac{ b_1 u_1-2u_2}{  b_1}\right) s\right),& 
				 \alpha\in \{\alpha_1,\alpha_1+2\alpha_2 \} \\
			\end{cases}
		\end{align}
		that is, in terms of $\fD=\fD_1\fD_2$, 
		  
		  \begin{align} 
			\textrm{Res}^2_x\left(p_b^\ast \frac{d_b \fD(\bu)}{\fD(b)} \ s\right) =
			\begin{cases}
				\textrm{Res}^2_x\left(\frac{1}{\bal^2}	p_b^\ast \left(\frac{(d\fD_1)_b(\bu)}{b_1}\right) s\right)
				& \textrm{ for } \alpha \textrm{ short}\\
				\frac{1}{2} \textrm{Res}^2_x\left(\frac{1}{\bal^2}	p_b^\ast \left(\frac{(d\fD_2)_b(\bu)}{ b_1}\right) s\right)
				 & \textrm{ for } 
				\alpha \textrm{ long}\\
			\end{cases}
		\end{align}

	  \end{Prop}

	  Again, as in the $G_2$ and $A_2$ cases, the result holds without requiring $s$ to be $W$-invariant.

		\begin{Prop}\label{nabla_lambda_B2}
		Consider
	   	 $b=(b_1,b_2)\in \scB$,   $\bv=(v_1,v_2)\in T_b\scB=\cB_{\fs\fo_5}$ and a ramification point 
	   	$x\in D_\alpha$ of $p_b:\widetilde{X}_b\to X$. Let
	   	$(U,\psi)$ be a local chart on $X$ and $z$ the respective
	   	   local coordinate,  centered at $z_0=p_b(x)$. Then

	   	\begin{align}
	   		(\psi^{-1})^\ast(\nabla^{GM}_\bv\bla_{SW})_b (x)=
			    \begin{cases}
			    	\frac{(d\fD_1)_b(\bv)[dz]}{j^1(\fD_1(b))}(z_0) [d\alpha]_x \ctimes
						    \begin{cases}
	   	 	      	\begin{pmatrix}
	   	 	      		-1\\ 1\\
	   	 	      	\end{pmatrix}, 
				      & \bal = \bal_2\\
	   	 	      	\begin{pmatrix}
	   	 	      		1\\ 0\\
	   	 	      	\end{pmatrix}, 
				     & \bal= \bal_1+\bal_2\\
						    \end{cases} \\
						    \\
				\frac{(d\fD_2)_b(\bv)[dz]}{j^1(\fD_2(b))}(z_0) [d\alpha]_x \ctimes
					      \begin{cases}
					      \frac{1}{2}
					         \begin{pmatrix}
	   	 	      		2\\ -1 \\
	   	 	      	\end{pmatrix}, & \bal =  \bal_1 \\
					\frac{1}{2}
					   \begin{pmatrix}
	   	 	      		0\\ 1\\
	   	 	      	\end{pmatrix}, & \bal=  \bal_1+2\bal_2\\
					      \end{cases}
			    \end{cases}
	   	\end{align}
	    and, explicitly, 
		    
	   	\begin{equation}
		 \frac{(d\fD_1)_b(\bv)[dz]}{j^1(\fD_1(b))}(z_0)= \frac{v_2 [dz]}{j^1(b_2)}(z_0),
	\end{equation}
	\begin{equation}
		 \frac{(d\fD_2)_b(\bv)[dz]}{j^1(\fD_2(b))}(z_0)= \frac{(2b_1 v_1-4b_2)[dz]}{j^1\left(b_1^2-4b_2 \right)}(z_0).
	\end{equation}

	   	 \end{Prop}

	   	 \emph{Proof of Theorem \ref{thm4}:}
	   	 
	   	 The Killing form $\kappa$ has matrix $6\begin{pmatrix}
	   	                               	1&1\\
	   	                               	1&2\\
	   	                               \end{pmatrix}$ with respect to the chosen basis. Consequently, 
	   	 $\begin{pmatrix}
	   	  	1\\ -\frac{1}{2}\\
	   	  \end{pmatrix}$ squares to $3$, $\begin{pmatrix}
						      0\\ 1\\
						  \end{pmatrix}$ squares to $12$ and both
		$\begin{pmatrix}
		 	1\\ -1\\
		 \end{pmatrix}$ and $\begin{pmatrix}
				      1\\ 0\\	  
				     \end{pmatrix}$ square to $6$ with respect to $\kappa$. Then Proposition \ref{nabla_lambda_B2}
		implies that for $x\in D_\alpha$ one has
		
		 \begin{multline}
		    (\psi^{-1})^\ast  (\nabla^{GM}_\bv\bla_{SW})_b\cup (\nabla^{GM}_\bw\bla_{SW})_b (x)= \\			    \begin{cases} \displaystyle
			    	6\frac{(d\fD_1)_b(\bv)(d\fD_1)_b(\bw)[dz]^2}{j^1(\fD_1(b))^2}(z_0) [d\alpha]^2_x,& \alpha \textrm{ short}\\[8pt] \displaystyle
			    	3 \frac{(d\fD_2)_b(\bv)(d\fD_2)_b(\bw)[dz]^2}{j^1(\fD_2(b))^2}(z_0) [d\alpha]^2_x,& \alpha \textrm{ long}.
			    \end{cases}
	   	 \end{multline}

	   	 Then Proposition \ref{residue_gen_B2} tells us that

	   	 \begin{multline}
	   	 	\textrm{Res}^2_x\left(p_b^\ast \frac{d_b \fD(\bu)}{\fD(b)} (\nabla^{GM}_\bv\bla_{SW})_b\cup (\nabla^{GM}_\bw\bla_{SW})_b  \right) =\\
	   	 	=\begin{cases} \displaystyle
			    	6\frac{(d\fD_1)_b(\bu)(d\fD_1)_b(\bv)(d\fD_1)_b(\bw)}{b_1 j^1(\fD_1(b))^2}(z_0),& \alpha \textrm{ short}\\[8pt] \displaystyle
			    	\frac{3}{2} \frac{(d\fD_2)_b(\bu)(d\fD_2)_b(\bv)(d\fD_2)_b(\bw)}{b_1 j^1(\fD_2(b))^2}(z_0),& \alpha \textrm{ long}\\
			    \end{cases} \\[8pt] \displaystyle
			 = \begin{cases} \displaystyle
			   	6 \frac{u_2 v_2 w_2}{b_1 j^1(b_2)^2}(z_0), & \alpha \textrm{ short}\\[8pt] \displaystyle
			   	12 \frac{(b_1 u_1-2u_2)(b_1 v_1- 2v_2)(b_1w_1- 2 w_2)}{b_1 j^1(b_1^2-4b_2)^2}(z_0),& \alpha \textrm{ long}.
			   \end{cases}.
	   	 \end{multline}
	   	 
	   	 The formula from the statement of  Theorem \ref{thm4} then follows, once we take into account that there are two points in the support of
	   	 each $D_\alpha$, that there are two short and two long simple roots, and that there is an overall factor of $\frac{1}{2}$ in
	   	 formula (\ref{BP_ver1}).
	   	 \qed

	 \end{document}